\documentclass[12pt]{article}
\usepackage{amsmath,amssymb}

\newcommand{\ncm}{\newcommand}
\ncm{\Inn}{\mbox{\rm Inn}} \ncm{\Ap}{\mbox{$\overline{\rm Inn}$}}
\ncm{\Sp}{\mbox{\rm Sp}} \ncm{\Ex}{\mbox{\rm Ex}}
\ncm{\OExt}{\mbox{\rm OrderExt}} \ncm{\AI}{\mbox{\rm AInn}}
\ncm{\HI}{\mbox{\rm HInn($A$)}} \ncm{\Aut}{\mbox{\rm Aut}}
\ncm{\Mal}{\mbox{$M_{\alpha}$}} \ncm{\Aff}{\mbox{${\rm Aff}$}}
\ncm{\id}{\mbox{\rm id}} \ncm{\Ker}{\mbox{\rm Ker}}
\ncm{\BE}{\begin{eqnarray*}} \ncm{\EE}{\end{eqnarray*}}
\ncm{\lra}{\mbox{$\longrightarrow$}} \ncm{\Hom}{\mbox{\rm Hom}}
\ncm{\calU}{{\cal U}} \ncm{\el}{\ell}
\ncm{\ad}{{\rm ad}}
\ncm{\diag}{{\rm diag}}
\ncm{\Ran}{{\rm Ran}}
\ncm{\Tr}{{\rm Tr}}
\ncm{\Arg}{{\rm Arg}}
\ncm{\Alg}{\mbox{\rm Alg}} \ncm{\Conv}{\mbox{\rm Conv}}
\ncm{\D}{{\mathcal {D}}} \ncm{\calE}{{\mathcal{E}} }    \ncm{\calJ}{{\mathcal{J}} }
\ncm{\calS}{{\mathcal{S}} } \ncm{\UHF}{{\rm UHF}} \ncm{\K}{{\mathcal
{K}}}\ncm{\B}{{\mathcal {B}}}

\ncm{\cstar}{C$^{*}$-algebra} \ncm{\cstars}{C$^{*}$-algebras}
\ncm{\cstarsub}{C$^{*}$-subalgebra} \ncm{\ra}{\mbox{$\rightarrow$}}
\ncm{\la}{\mbox{$\leftarrow$}} \ncm{\hra}{\hookrightarrow}
\ncm{\da}{\mbox{$\downarrow$}} \ncm{\se}{\mbox{$\searrow$}}
\ncm{\al}{\mbox{$\alpha $}} \ncm{\del}{\mbox{$\delta$}}
\ncm{\supp}{\mbox{\rm supp}} \ncm{\Ad}{\mbox{\rm Ad}}
\ncm{\CAR}{\mbox{$M_{2^{\infty}}$}} \ncm{\ep}{\mbox{$\eps > 0$}}
\ncm{\ol}{\overline} \ncm{\Mninf}{\mbox{$M_{n^{\infty}}$}}
\ncm{\MR}{M. R\o{}rdam} \ncm{\Range}{\mbox{\rm Range}}
\ncm{\vo}{}%{\bf}
\ncm{\ch}{}%{\it}
\ncm{\CMP}{Comm. Math. Phys.} \ncm{\add}{}
\ncm{\tilsig}{\tilde{\sigma}} \ncm{\dist}{{\rm
dist}}\ncm{\eps}{\varepsilon}  \ncm{\calL}{{\mathcal{L}}}
\ncm{\E}{{\mathcal{E}} }    \ncm{\M}{{\mathcal{M}} }
\ncm{\calO}{{\mathcal{O}}} \ncm{\F}{{\mathcal{F}}} \ncm{\G}{{\mathcal{G}} }
\ncm{\calH}{{\mathcal{H}}}   \ncm{\calC}{{\mathcal{C}}}
\ncm{\lan}{{\langle}}\ncm{\ran}{{\rangle}}
\ncm{\calK}{{\mathcal{K}}}    \ncm{\Spec}{{\rm Spec}}
\ncm{\calP}{{\mathcal{P}}} \ncm{\Hil}{{\mathcal{H}}}
\ncm{\U}{{\mathcal{U}}}                \ncm{\A}{\mathcal{A}}
\ncm{\Nr}{\mathcal{N}}  %\ncm{\calC}{\mathcal{C}}

\newtheorem{theo}{Theorem}[section]
\newtheorem{cor}[theo]{Corollary}
\newtheorem{lem}[theo]{Lemma}
\newtheorem{prop}[theo]{Proposition}
\newtheorem{remark}[theo]{Remark}
\newtheorem{definition}[theo]{Definition}
\newtheorem{example}[theo]{Example}
\newtheorem{property}[theo]{Property}

\newenvironment{pf}{{\it Proof.}}{\hfill$\square$\vspace{3mm}}
\newenvironment{pff}{{\it Proof}}{\hfill$\square$\vspace{3mm}}

\ncm{\R}{\mbox{\bf R}} \ncm{\Z}{\mbox{\bf Z}} \ncm{\T}{{\bf T}}
\ncm{\TT}{\T$^{2}$} \ncm{\N}{\mbox{\bf N}} \ncm{\C}{\mbox{\bf C}}
\ncm{\Dk}{{\bf D}}

\title{Quasi-diagonal flows}
\bigskip

\author{ Akitaka Kishimoto \\ {\small
Department of Mathematics, Hokkaido University, Sapporo 060-0810,
Japan}\\ and\\
Derek W. Robinson\\
{\small Centre for Mathematics and its Applications, Australian
National University}\\ {\small Canberra, ACT 0200, Australia}}

\date{December 2008}

\begin{document}
\maketitle

\begin{abstract}
We introduce two notions for flows on quasi-diagonal \cstars,
quasi-diagonal and pseudo-diagonal flows; the former being
apparently stronger than the latter. We derive basic facts about
these flows and give various examples. In addition we  extend
results of   Voiculescu from  quasi-diagonal \cstars\ to these
flows.
\medskip


MSC(2000) 46L55
\end{abstract}

\section{Introduction}
Flows on \cstars\ have been studied for some time; basic facts on
flows and their generators, from the perspectives of functional
analysis, spectral analysis, and Hilbert space representation
theory, etc. are described in \cite{BR1,BR2}. But there remain many
problems pertaining to \cstars. For example we still lack  clear and
useful criteria which  distinguish various kinds of flows, e.g.
approximately inner flows and,  in the case of AF algebras,
(approximate) AF flows. (See \cite{Sak} for some results for flows
on AF algebras.) We hope to contribute towards clarification of the
situation  by introducing other properties of flows which appear  to
have close bearing on these features at least in the case of simple
\cstars.

A bounded operator $T$ on a separable Hilbert space $\Hil$
is called {\em quasi-diagonal} if there is an increasing
sequence $(E_n)$ of finite-rank projections on $\Hil$ such that
$\lim_n E_n=1$ strongly and $\|[E_n,T]\|\to 0$.
This notion is extended to a norm-closed $^*$-algebra $A$ of bounded operators:
$A$ is called quasi-diagonal if there is such a sequence $(E_n)$ and
$\|[E_n,T]\|\to 0$ for all $T\in A$.
If $A$ is a C$^*$-algebra,
then $A$ is called quasi-diagonal if there is a faithful
representation $\pi$ of $A$ such that $\pi(A)$ is quasi-diagonal.
(See \cite{Voi93,Br04} for more details.)
Easy examples of
quasi-diagonal \cstars\ include AF algebras and commutative \cstars.
We mimic this notion in application to flows on \cstars\ in two ways.

\begin{definition}
Given a Hilbert space $\Hil$, let $A$ be a norm-closed $^*$-algebra
of bounded operators on $\Hil$ and $U$ a unitary flow on
$\Hil$ such that $U_txU_t^*\in A$ for $t\in\R$ and $t\mapsto
U_txU_t^*$ is norm-continuous for any $x\in A$.

We call $(A,U)$  {\em quasi-diagonal} if for any finite set
$\F$ of $A$, any finite set $\omega$ of $\Hil$ and $\eps>0$
there is a finite-rank projection $E$ on $\Hil$ such that
$\|[E,x]\|\leq\eps\|x\|$ for $x\in \F$,
$\|(1-E)\xi\|\leq\eps\|\xi\|$ for $\xi\in \omega$  and
$\|[E,U_t]\|<\eps$ for $t\in [-1,1]$.
We call $(A,U)$  {\em
pseudo-diagonal} if for any finite set $\F$ of $A$, any finite set
$\omega$ of $\Hil$, and $\eps>0$ there is a finite-rank
projection $E$ on $\Hil$ and a unitary flow $V$ on $E\Hil$ such that
$\|[E,x]\|\leq\eps\|x\|$ for $x\in \F$,
$\|(1-E)\xi\|\leq\eps\|\xi\|$ for $\xi\in \omega$  and
$\|EU_txU_t^*E-V_tExEV_t^*\|\leq\eps\|x\|$ for $x\in\F$ and $t\in
[-1,1]$.

Let $A$ be a \cstar\ and let $\alpha$ be a flow on $A$. We call
$\alpha$  {\em quasi-diagonal} $($resp.\  {\em pseudo-diagonal}$)$ if
$(A,\alpha)$ has a covariant representation $(\pi,U)$ on a Hilbert
space $\Hil_\pi$, with $\pi$ faithful and non-degenerate, such that
$(\pi(A),U)$ is quasi-diagonal $($resp.\  pseudo-diagonal$)$.
\end{definition}

In the above definition $\pi$  is required  to be non-degenerate.
But this is not essential.
A direct proof will be given in the
beginning of Section 2 but this also follows from Theorems~\ref{Ch}
and \ref{Ch'} below.
Thus we immediately obtain the following result.
(We  do not know if a similar statement is true or false for  approximately
inner flows.)

\begin{cor}
Let $\alpha$ be a quasi-diagonal $($resp.\  pseudo-diagonal$)$ flow on a
\cstar\ $A$ and $B$ an $\alpha$-invariant \cstarsub\ of $A$.
Then $\alpha|_B$ is quasi-diagonal $($resp.  pseudo-diagonal$)$.
\end{cor}

Let $H$ denote the self-adjoint generator of $U$ in the above
definition.
In general $H$ is unbounded.
If $Q$ is a a bounded operator on $\Hil$ then  $[Q,H]$ is
defined to be bounded if $Q\D(H)\subset
\D(H)$ and $QH-HQ$ is bounded on $\D(H)$ (and extends to a bounded
operator on $\Hil$).
We may replace the condition $\|[E,U_t]\|<\eps$
for $t\in [-1,1]$ in the definition of  quasi-diagonality by
the seemingly stronger condition $\|[E,H]\|<\eps$.
The opposite implication can be seen
from the proposition given below.
Using this we conclude that
quasi-diagonality implies  pseudo-diagonality since if
$\|[E,H]\|<\eps$ and we set $V_t=e^{itEHE}$ then
$\|EU_t\pi(x)U_t^*-V_tE\pi(x)EV_t^*\|\leq 2\eps\|\pi(x)\|$ for
any $x\in A$.

\begin{prop}\label{smear}
Let $U$ be a flow on $\Hil$ and  $H$  the self-adjoint
generator of $U$.
For any $\eps>0$ there is a $\delta>0$
satisfying the following condition.

If $E$ is a projection such that
$\|[E,U_t]\|<\delta$ for $t\in [-1,1]$, then there is a projection
$F$ on $\Hil$ such that $\|E-F\|<\eps$ and $\|[F,H]\|<\eps$.
\end{prop}
\begin{pf}
Note that it follows from the above estimate on $\|[E,U_t]\|$ that
$\|[E,U_t]\|<\delta(1+|t|)$ for all $t\in\R$.
In addition to  this estimate
we  use the fact that $t\mapsto U_tEU_t^*$ is continuous in the
strong operator topology.

Let $f$ be a non-negative $C^\infty$ function on $\R$ such that
$\supp(f)\subset [1/3,4/3]$ and $f(t)=1$ for $t\in [2/3,1]$. Define
$\hat{f}$ by $\hat{f}(p)=(2\pi)^{-1}\int e^{-ipt}f(t)dt$ and set
$C=\int|t\hat{f}(t)|dt<\infty$.
Let $g$ be a non-negative $C^\infty$
function on $\R$ such that the support of $g$ is compact, $\int
g(t)dt=1$ and $\int |g'(t)|dt<\eps/C$.
Set $D=\int g(t)(1+|t|)dt$.
Assuming $\delta D<\eps/2<1/3$  we define
$$
Q=\int g(t)U_tEU_t^*dt\;.
$$
Then $0\leq Q\leq 1$, $\|Q-E\|<\eps/2$  and
$\|[H,Q]\|<\eps/C$, where $i[H,Q]$ is identified with $-\int
g'(t)U_tEU_t^*dt$.
Since $\Sp(Q)\subset [0,\eps/2)\cup
(1-\eps/2,1]$ it follows that $F=f(Q)=\int \hat{f}(t)e^{itQ}dt$ is
a projection satisfying $\|F-Q\|<\eps/2$.
It also follows
that $\|[H,F]\|\leq \|[H,Q]\|\int |t\hat{f}(t)|dt<\eps$.
Since
$\|E-F\|<\eps$, this concludes the proof. (See \cite{BR1} for
the norm estimates used here.)
\end{pf}

We note that a covariant representation $(\rho,V)$ of $(A,\alpha)$
naturally induces a representation $\rho\times V$ of the crossed
product $A\times_\alpha\R$ on the representation space $\Hil_\rho$
of $\rho$. We denote by $\K(\Hil_\rho)$ the compact operators on
$\Hil_\rho$.

By extending Voiculescu's theorem \cite{Voi76} to accommodate the
flow we establish the following:

\begin{theo}\label{Voi}
Let $\alpha$ be a quasi-diagonal $($resp.\  pseudo-diagonal$)$ flow on
$A$.
If $(\rho,V)$ is a covariant representation
of $A$ such
that $\rho\times V$ is a faithful representation of
$A\times_\alpha\R$ and $\Ran(\rho\times V)\cap \K(\Hil_\rho)=\{0\}$ then
$(\rho(A),V)$ is quasi-diagonal $($resp.\ pseudo-diagonal$)$.
\end{theo}

Mimicking the corresponding result due to Voiculescu \cite{Voi91} we
shall give characterizations of quasi-diagonal and pseudo-diagonal
flows.

If  $A$ and $B$ are \cstars\  then a linear map $\phi$ of $A$ into
$B$ is called {\em positive} if $\phi(A_+)\subset B_+$ and {\em completely positive} (or CP)
if $\phi_n=\id\otimes \phi\colon M_n\otimes A \to M_n\otimes B$ is positive for all $n$.

\begin{theo}\label{Ch}
Let $\alpha$ be a flow on a \cstar\ $A$. Then the following
conditions are equivalent:
\begin{enumerate}
\item $\alpha$ is quasi-diagonal.
\item For any finite subset $\F$ of $A$ and $\eps>0$ there is
a finite-dimensional \cstar\ $B$,  a flow $\beta$ on $B$  and a CP
map $\phi$ of $A$ into $B$ such that $\|\phi\|\leq1$,
$\|\phi(x)\|\geq (1-\eps) \|x\|$ and
$\|\phi(x)\phi(y)-\phi(xy)\|\leq\eps\|x\|\|y\|$ for $x,y\in\F$,
and $\|\beta_t\phi-\phi\alpha_t\|<\eps$ for $t\in[-1,1]$.
\item For any finite subset $\F$ of $A$ and $\eps>0$ there is a
covariant representation $(\pi,U)$ and  a finite-rank
projection $E$ on $\Hil_\pi$ such that $\|E\pi(x)E\|\geq
\|x\|-\eps$ and $\|[E,\pi(x)]\|\leq\eps\|x\|$ for $x\in\F$
and $\|[E,U_t]\|<\eps$ for $t\in [-1,1]$.
\end{enumerate}
\end{theo}

\begin{theo}\label{Ch'}
Let $\alpha$ be a flow on a \cstar\ $A$. Then the following
conditions are equivalent:
\begin{enumerate}
\item $\alpha$ is pseudo-diagonal.
\item For any finite subset $\F$ of $A$ and $\eps>0$ there is
a finite-dimensional \cstar\ $B$, a flow $\beta$ on $B$ and a CP
map $\phi$ of $A$ into $B$ such that $\|\phi\|\leq1$,
$\|\phi(x)\|\geq (1-\eps) \|x\|$ and
$\|\phi(x)\phi(y)-\phi(xy)\|\leq\eps\|x\|\|y\|$ for $x,y\in\F$,
and $\|\beta_t\phi(x)-\phi\alpha_t(x)\|\leq\eps\|x\|$ for
$x\in\F$ and $t\in[-1,1]$.
\item For any finite subset $\F$ of $A$ and $\eps>0$ there is a
covariant representation $(\pi,U)$, a finite-rank projection $E$ on
$\Hil_\pi$ and a unitary flow $V$ on $E\Hil_\pi$ such that
$\|E\pi(x)E\|\geq (1-\eps)\|x\|$ and
$\|[E,\pi(x)]\|\leq\eps\|x\|$ for $x\in\F$ and
$\|EU_t\pi(x)U_t^*E-V_tE\pi(x)EV_t^*\|\leq\eps\|x\|$ for $x\in
\F$ and $t\in [-1,1]$.
\end{enumerate}
\end{theo}

In the above theorems the finite-dimensional \cstar\ $B$ can be
assumed to be a matrix algebra $M_k$ for some $k\in \N$.

If  $A$ is separable and $\alpha$ is a pseudo-diagonal flow let
$(\F_n)$ be an increasing sequence of finite subsets of $A$ whose
union is dense in $A$ and choose, for each $(\F_n,n^{-1})$ in place
of $(\F,\eps)$, a CP map $\phi_n$ into $M_{k_n}$ and a flow
$\beta^{(n)}$ on $M_{k_n}$ as specified in condition (2) of the
above theorem. Thus we can define a non-continuous flow $\beta$ on
the direct product $B=\prod_nM_{k_n}$ by
$\beta_t(x)=\prod_n\beta^{(n)}_t(x_n)$ for $x=(x_n)\in B$ and a CP
map $\phi$ of $A$ into $B$ by $\phi(x)=(\phi_n(x))_n$. Let
$I=\bigoplus_nM_{k_n}$, which is the ideal of $B$ consisting of
sequences converging to zero, and let $Q$ denote the quotient map of
$B$ onto $B/I$. Then it follows that $\psi=Q\phi$ is an isomorphism
of $A$ into $B/I$  satisfying $\psi\alpha_t=\beta_t\psi$. A
separable \cstar\ is an MF algebra if it can be embedded into
$\prod_n M_{k_n}/\bigoplus_nM_{k_n}$ for some $(k_n)$ (see \cite{BK}
for MF algebras). We may call the flow $\alpha$ an {\em MF flow}
since it satisfies the intertwining property with $\beta$. It might
be interesting to explore this class of flows.

We will show that if $\alpha$ is an approximately inner flow on a
quasi-diagonal \cstar\
then $\alpha$ is pseudo-diagonal (Proposition~\ref{AI}).
We will also show that if $\alpha$ is a pseudo-diagonal
flow on a unital \cstar\
then $\alpha$ has KMS states for all
inverse temperatures (Proposition~\ref{KMS}).

If $A$ is an AF algebra and $\alpha$ is an (approximate) AF flow
then it follows that $\alpha$ is quasi-diagonal (Proposition~\ref{AF}). If
$A$ is an AF algebra which has a faithful family of type I
quotients then any flow on $A$ is quasi-diagonal (Proposition~ \ref{TI}).

Let $\alpha$ (resp.\  $\beta$) be a flow on a \cstar\ $A$ (resp.\  $B$).
We say that $(B,\beta)$ homotopically dominates $(A,\alpha)$ if
there are homomorphisms $\phi\colon A \to B$ and $\psi\colon B \to A$ and a
homotopy $\{\chi_s;\ s\in [0,1]\}$ of homomorphisms of $A$ into $A$
such that $\phi\alpha_t=\beta_t\phi$, $\psi\beta_t=\alpha_t\psi$,
$\chi_s\alpha_t=\alpha_t\chi_s$, $\chi_0=\psi\phi$ and
$\chi_1=\id_A$.
The main result of Voiculescu's paper \cite{Voi91}
has the following analogue:

\begin{theo}\label{hom}
Suppose that $(B,\beta)$ homotopically dominates $(A,\alpha)$.
If $\beta$ is quasi-diagonal then $\alpha$ is
quasi-diagonal.
\end{theo}

This implies that if $\alpha$ is a flow on a \cstar\ $A$ then the flow
$\alpha\otimes\id$ on $A\otimes C_0[0,1)$ is quasi-diagonal.
(The family of endomorphisms $\phi_s,s\in[0,1]$ of $A\otimes C_0[0,1)$
defined by $\phi_s(x)(t)=x(st)$ commutes with the flow
$\alpha\otimes\id$ and satisfies $\phi_1=\id$ and $\phi_0=0$.
This also follows directly from Proposition \ref{hom1}.)
Thus approximate innerness does not follow from quasi-diagonality without
further conditions on the \cstar.
Another result of this type is that if
$\alpha$ is a flow on a quasi-diagonal \cstar\ $A$  then the flow
$\beta$ on $A\otimes C[0,1]$ defined by
$\beta_t(x)(s)=\alpha_{st}(x(s))$ is quasi-diagonal (Proposition~\ref{hom2}).

We note that we have not been able to give the pseudo-diagonal
version of the above theorem. We also note that we do not know if
quasi-diagonality is strictly stronger than pseudo-diagonality or
not.

Let $u$ be an $\alpha$-cocycle, i.e.\  let $u$ denote a continuous
function from $\R$ into the unitary group of $M(A)$ such that
$t\mapsto u_t$ is continuous in the strict topology and
$u_s\alpha_s(u_t)=u_{s+t}$ for $s,t\in\R$.
If $A$ is unital then
the multiplier algebra $M(A)$ is just $A$ and the strict topology is
the norm topology.
We say the flow $t\mapsto \Ad\,u_t\alpha_t$ is a
cocycle perturbation of $\alpha$.
We note that  quasi-diagonality
(resp.\  pseudo-diagonality)  is stable under cocycle perturbations
(Propositions~\ref{pert} and  \ref{PS}).
We also note that if $B$ is
an $\alpha$-invariant hereditary C$^*$-subalgebra of $A$ which
generates $A$ as an ideal  then $\alpha|_B$ is quasi-diagonal (resp.\
pseudo-diagonal) if and only if $\alpha$ is quasi-diagonal (resp.\
pseudo-diagonal) (Corollary~\ref{Her}).

In Section 2 we will give the above basic facts on quasi-diagonal
and pseudo-diagonal flows and some examples including the proof of
Theorem~\ref{hom}.
For example the rotation flow on the
continuous functions on the unit circle is not quasi-diagonal (and
not even pseudo-diagonal) but the rotation flow on the continuous
functions on the unit disk is quasi-diagonal.
In Section~3 we generalize Voiculescu's Weyl-von Neumann theorem \cite{Voi76} to
cover the present situation and thereby prove Theorem~\ref{Voi}.
In Section~4 we deal with the adaptation of Voiculescu's results in
\cite{Voi91} to prove Theorems \ref{Ch} and \ref{Ch'}.

\section{Quasi-diagonal and pseudo-diagonal flows}

Let $\alpha$ be a flow on $A$.
The definition of quasi-diagonality, or pseudo-diagonality, of  $\alpha$ required
the representation $\pi$ in the covariant representation $(\pi,U)$ of $(A,\alpha)$
to be  faithful and  non-degenerate.
But  the non-degeneracy of $\pi$ is not essential
by the following argument.

First this is evident if $A$ is unital.
Therefore we assume that $A$ is not unital.

Secondly, let $\pi$ be a faithful degenerate representation of
$A$ on a Hilbert space $\Hil$ and $U$ a unitary flow  on $\Hil$ such
that $\Ad\,U_t\pi(x)=\pi\alpha_t(x)$  for all $ x\in A$.
Let $P$ be the projection onto the closure of $\pi(A)\Hil$.
Note that $U_tP=PU_t$ and let us denote by $UP$ the unitary flow $t\mapsto U_tP$ on
$P\Hil$.

Suppose that $(\pi(A),U)$ is pseudo-diagonal.
We shall show that
$(\pi(A)P, UP)$ is pseudo-diagonal.
For a finite subset $\F$ of $A$,
a finite subset $\omega$ of $P\Hil$ and $\eps>0$  we choose a
finite-rank projection $E$ on $\Hil$ and a unitary flow $V$ on
$E\Hil$ which satisfy the conditions of the definition.
Let $\K_1$ be the subspace $(1-P)E\Hil$.
We find a subspace $\K_2$ of
$P\Hil$ with the same dimension as $\K_1$ such that $\K_2$ is
orthogonal to $PE\Hil$ and
$\|\pi\alpha_t(x)|_{\K_2}\|\leq(\eps/2)\|x\|$ for $x\in \F\cup\F^*$
and $t\in [-1,1]$.
Let $W_1$ be a unitary from $\K_1$ onto $\K_2$
and denote by $P_i$ the projection onto $\K_i$ for $i=1,2$.
Regarding $W_1$ as $W_1=W_1P$, let $W=W_1+W_1^*+(1-P_1-P_2)$,  which
is a unitary on $\Hil$, and let $F=WEW^*$. Since $WE\Hil\subset
W(1-P)E\Hil+WPE\Hil\subset P_2W(1-P)E\Hil+PE\Hil$, it follows that
$F\leq P$.
Since $\pi(x)WE=\pi(x)W_1E+\pi(x)(1-P_2)E=\pi(x)P_2(W_1-1)E+\pi(x)E$  we
obtain  $\|\pi(x)WE-\pi(x)E\|<\eps$, which implies that
$\|F\pi\alpha_t(x)F-WV_tW^*F\pi(x)FWV_t^*W^*\|\leq 5\eps\|x\|$
for $x\in\F$ and $t\in [-1,1]$.
The other properties follow easily.
Thus the pair $F$ and $t\mapsto WV_tW^*$ satisfies the required
conditions for $(\pi(A)P,UP)$.

Now suppose that $(\pi(A),U)$ is quasi-diagonal. Let
$(\bar{\pi},\bar{U})$ be the direct sum of $(\pi,\chi_pU)$ over all
rational numbers $p$, on the representation space
$\bar{\Hil}=\bigoplus_p\Hil$, where $\chi_pU$ is the unitary flow
$t\mapsto e^{ipt}U_t$. Let $P$ be the projection onto the closure of
$\pi(A)\Hil$ as before and let $\bar{P}$ be the projection onto the
closure of $\bar{\pi}(A)\bar{\Hil}$, i.e.\  $\bar{P}=\bigoplus_pP$.
We shall show that $(\bar{\pi}(A)\bar{P},\bar{U}\bar{P})$ is
quasi-diagonal.

From now on we use $\pi, U,P$ to  denote
$\bar{\pi},\bar{U},\bar{P}$. We have now assumed that $\pi\times U$
is faithful besides $(\pi(A),U)$ being quasi-diagonal. Let $H$ be
the self-adjoint generator of $U$. For a finite subset $\F$ of $A$,
a finite subset $\omega$ of $P\Hil$ and $\eps>0$ we choose a
finite-rank projection $E$ on $\Hil$ such that $\|[E,H]\|<\eps$
holds in addition to  the other conditions in the definition. There
is a finite-dimensional subspace $\K_1$ of $(1-P)\Hil$ such that
$\K_1\supset (1-P)E\Hil$ and $\|[P_1,H]\|<\eps/2$, where $P_1$ is
the projection onto $\K_1$. Let
$\lambda_1,\lambda_2,\ldots,\lambda_n$ be the eigenvalues of
$P_1HP_1$ in increasing order. We choose a finite-rank projection
$P_2$ such that $P_2\leq P$, $P_2PE=0$,
$\|P_2\pi(x)\|,\|\pi(x)P_2\|\leq (\eps/2)\|x\|$ for $x\in\F$,
$\|[P_2,H]\|<\eps/2$ and the increasing list of eigenvalues of
$P_2HP_2$ are arbitrarily close to
$\lambda_1,\lambda_2,\ldots,\lambda_n$. (In particular $P_1$ and
$P_2$ have the same rank.) This is possible by the lemma below which
uses faithfulness of $\pi\times U$. Then we choose a unitary $W_1$
of $\K_1$ onto $P_2\Hil$ such that $W_1P_1HP_1\approx P_2HP_2W_1$.
We set $W=W_1+W_1^*+(1-P_1-P_2)$. Then $F=WEW^*\leq P$ and
$\|U_tW-WU_t\|\leq \eps |t|$ (by making $W_1P_1HP_1\approx
P_2HP_2W_1$ precise). This implies that $F$ satisfies the required
conditions.

\begin{lem}
Suppose that $A$ is non-unital and let $\pi,U,P$ be as above.  For
any finite subset $\F$ of $A$, $\lambda\in\R$  and $\eps>0$ there
exists a unit vector $\xi\in P\Hil$ such that $\|\pi(x)\xi\|<\eps$
for $x\in\F$ and $\|U_t\xi-e^{i\lambda t}\xi\|\leq \eps|t|$.
\end{lem}
\begin{pf}
Let $z=\sum_{x\in\F}x^*x$ and let $P_H$ denote the spectral measure
for $H$. Suppose that there is an $\eps>0$ such that $\lan
\xi,\pi(z)\xi\ran\geq \eps$ for any unit vector $\xi$ in
$P_H(\lambda-\eps,\lambda+\eps)P\Hil$. Let $\hat{f}$ be a
non-negative $C^\infty$-function on $\R$ such that $f\not=0$ and
$\supp(\hat{f})\subset (\lambda-\eps,\lambda+\eps)$. Since
$\pi\times U$ is faithful $ \lambda(f)(z-\eps)\lambda(f)^*\geq0, $
where $\lambda(f)=\int f(t)\lambda_t dt$ is a multiplier of
$A\times_\alpha\R$ such that $\pi(\lambda(f))=\hat{f}(H)$. Applying
$\hat{\alpha}_p$ and taking the integral over $p$ implies that
$\int|f(t)|^2 \alpha_t(z)$ is invertible, which contradicts that $A$
is non-unital. (See 7.8 of \cite{Ped} for more details.)
\end{pf}

In order for $\alpha$ to be quasi-diagonal or pseudo-diagonal the
\cstar\ $A$ must be quasi-diagonal.
Moreover,  it follows that if
$\alpha$ is quasi-diagonal then the crossed product
$A\times_\alpha\R$ is quasi-diagonal.
(The pair $(\pi,U)$ gives a
representation $\pi\times U$ of $A\times_\alpha\R$, which may not be
faithful, such that $\pi\times U(A\times_\alpha\R)$ is
quasi-diagonal.
As a faithful representation of $A\times_\alpha\R$ is
required in the definition of quasi-diagonality we may take the
direct sum of $\pi\times \chi_p U$ over all rationals $p$ as in the
previous paragraph.)

If  $\alpha$ is a flow on an AF algebra $A$ then the crossed
product $A\times_\alpha\R$ is AF-embeddable; in particular it is
quasi-diagonal.
(We learned this fact from M. Izumi; the argument
uses the fact that the crossed product of $A$ by $\alpha|_{\Z}$ is
AF-embeddable, due to \cite{Voi86} and \cite{Br98}.)

As we shall see $\alpha$ need not be quasi-diagonal, nor
pseudo-diagonal,  even if $A\times_\alpha\R$ is quasi-diagonal.

\begin{prop}\label{pert}
Let $\alpha$ be a flow on $A$ and let $u$ be an $\alpha$-cocycle.
Then $\alpha$ is quasi-diagonal if and only if $t\mapsto
\Ad\,u_t\alpha_t$ is quasi-diagonal.
\end{prop}
\begin{pf}
If $A$ is unital this follows straightforwardly. Suppose that $A$
does not have a unit and that $\alpha$ is quasi-diagonal. Thus we
assume that $A$ acts on a Hilbert space $\Hil$ non-degenerately and
there is a unitary flow $U$ such that $\alpha_t(x)=U_t xU_t^*$ for
$x\in A$ and $(A,U)$ is quasi-diagonal. Let $\F$ be a finite subset
of $A$ and $\omega$ a finite subset of $\Hil$. Then we choose
$p,e\in A$ such that $0\leq p\leq e\leq 1$, $ep=p$,
$\|x-pxp\|\approx0$ for $x\in \F$, $\|p\xi-\xi\|\approx0$ for
$\xi\in\omega$  and $\|\alpha_t(e)-e\|\approx0$ for $t\in[-1,1]$. We
choose an $\alpha$-cocycle $v$ in $A+\C1$ such that
$\|(u_t-v_t)e\|\approx0,\ t\in [-1,1]$, where $t\mapsto v_t$ is
continuous in norm \cite{K06}. We choose a finite-rank projection
$E$ such that $\|[E,x]\|\approx0$ for $x\in\F\cup \{p,e\}\cup
\{v_t;\ t\in [-1,1]\}$, $\|(1-E)\xi\|\approx0$ for $\xi\in \Omega$
and $\|[E,H]\|\approx0$. By the lemma below there is a subprojection
$F$ of $E$ such that $Fp\approx Ep$, $Fe\approx F$, and
$\|[F,H]\|\approx 0$. Since $\|x-pxp\|\approx0$ for $x\in\F$  we
have  $\|[F,x]\|\approx \|[E,x]\|\approx0$ for $x\in\F$. Since
$u_te\approx v_te$ and $eu_t\approx ev_t$ we have
$\|[F,u_t]\|\approx \|[F,v_t]\|\approx0$ for $t\in [-1,1]$, which
implies that $\|[F,u_tU_t]\|\approx0$ for $t\in [-1,1]$. Further we
have $\|(1-F)\xi\|\approx \|(1-E)\xi\|\approx0$ for $\xi\in \omega$.
\end{pf}

\begin{lem}\label{pert1}
For any $\eps>0$ there exists a $\delta>0$ such that the following holds.

If $e,p\in A$ and a finite-rank projection $E$ satisfy $0\leq
p\leq e\leq 1$, $ep=p$, $\|\alpha_t(e)-e\|<\delta$ for $t\in
[-1,1]$, $\|[E,H]\|<\delta$, $\|[E,e]\|<\delta$ and
$\|[E,p]\|<\delta$, then there is a finite rank projection $F$ such
that $F\leq E$, $\|Ep-Fp\|<\eps$, $\|Fe-F\|<\eps$ and
$\|[F,H]\|<\eps$.
\end{lem}
\begin{pf}
Let $e'=EeE$, $p'=EpE$, and $H'=EHE$.
Since $\|[e^{itH},e']\|\approx0$ for $t\in [-1,1]$ and $\|(1-E)HE\|\approx0$
we conclude that $\|[e^{itH'},e']\|\approx 0$ for $t\in [-1,1]$.
Then Lin's theorem \cite{L97,BK01} for almost commuting self-adjoint
$e'$ and $H'$ in $\B(E\Hil)$ tells us that there is a self-adjoint
$h$ in $\B(E\Hil)$ such that $h\approx e'$ and $\|[q,H']\|\approx0$
uniformly for any spectral projection $q$ of $h$.
Since $p'e'\approx p'$we deduce that $p'h\approx p'$.
Let $F$ be the spectral projection of $h$ corresponding to $[1-\eps/2,1]$.
From the lemma below and $p'=EpE\approx Ep$ we may suppose that
$\|Fp-Ep\|<\eps$.
Since $\|Fh-F\|\leq \eps/2$ and $h\approx EeE$ we may also suppose
that $\|Fe-F\|<\eps$.
Since $[F,H]=[F,H']+FH(1-E)-(1-E)HF$  it follows that $\|[F,H]\|\leq
\|[F,H']\|+\|(1-E)HE\|$, which we may suppose is smaller than
$\eps$.
\end{pf}

\begin{lem}\label{proj}
For any $\eps,\eps'>0$ there is a $C>0$  such that the following holds.

For any $h,p\in A_{sa}$ such that $0\leq h\leq 1$, $0\leq p\leq 1$ and
$\|hp-p\|<\delta$ the spectral projection $F$ of $h$ corresponding
to $[1-\eps,1]$ satisfies  $\|Fp-p\|<\eps'+C\delta$.
\end{lem}
\begin{pf}
Fix a continuous function $f$ on $[0,1]$ such that $0\leq f\leq
1$, $f(1)=1$ and $\supp(f)\subset [1-\eps,1]$ and choose  a
polynomial $q(t)=\sum_{k=1}^n c_kt^k$ with $q(1)=1$ such that
$|f(t)-q(t)|<\eps'$ for $ t\in [0,1]$.
Since $\|q(h)p-p\|\leq \sum_{k=1}^n|c_k|k\delta\equiv C\delta$ and
$\|f(h)-q(h)\|<\eps'$ it follows that $\|(1-F)p\|=\|(1-F)(1-f(h))p\|<\eps'+C\delta$,
where $F$ denotes the spectral projection of $h$ corresponding to $[1-\eps,1]$.
\end{pf}

\begin{prop}\label{PS}
Let $\alpha$ be a flow on $A$ and $u$  an $\alpha$-cocycle.
Then $\alpha$ is pseudo-diagonal if and only if $t\mapsto
\Ad\,u_t\alpha_t$ is pseudo-diagonal.
\end{prop}
\begin{pf}
If $A$ is unital this follows straightforwardly. Suppose that $A$
does not have a  unit and that $\alpha$ is pseudo-diagonal. Thus we
assume that $A$ acts on a Hilbert space $\Hil$ non-degenerately and
there is a unitary flow $U$ such that $\alpha_t(x)=U_t xU_t^*,\ x\in
A$ and $(A,U)$ is pseudo-diagonal. Let $u$ be an $\alpha$-cocycle in
$M(A)$. Further let $\F$ be a finite subset of $A$ and $\omega$ a
finite subset of $\Hil$. Then, by the lemma below, we choose $p,e\in
A$ such that $0\leq p\leq e\leq 1$, $ep=p$, $\|x-pxp\|\approx0$ for
$x\in \F$, $\|p\xi-\xi\|\approx0$ for $\xi\in\omega$,
$\|\alpha_t(e)-e\|\approx0$ for $t\in[-1,1]$  and
$\|[e,u_t]\|\approx0$ for $t\in [-1,1]$. For $u$ and $e$ we choose
an $\alpha$-cocycle $v$ in $A+\C1$ such that
$\|(u_t-v_t)e\|\approx0$ for $t\in [-1,1]$ \cite{K06}. We then
express $v_t$ as $wU^{(h,\alpha)}_t\alpha_t(w^*)$, where $w\in
\U(A+\C1)$, $h=h^*\in A+\C1$ and $U_t^{(h,\alpha)}$ denotes the
$\alpha$-cocycle defined by ${d\over dt}
U^{(h,\alpha)}_t=U_t^{(h,\alpha)}\alpha_t(ih)$:
$$
U^{(h,\alpha)}_t=1+\sum_{n=1}^\infty
\int_{\Omega_n(t)}\alpha_{t_1}(ih)\alpha_{t_2}(ih)\cdots
\alpha_{t_n }(ih)dt_1\cdots dt_n
$$
where $\Omega_n(t)=\{(t_1,t_2,\ldots,t_n):\ 0\leq
t_1\leq t_2\leq\cdots\leq t_n\leq t\}$ for $t\geq 0$ and similarly
for $t<0$ (see Lemma 1.1 of \cite{K00}).

Let $\G=\F\cup \{p,e\}\cup \{\alpha_t(h),\alpha_t(w);\ t\in
[-1,1]\}$.
Then $\G$  is a compact subset of $A+\C1$.
Since $(A,U)$ is pseudo-diagonal, we choose, for $\G$ and $\omega$,
a finite-rank projection $E$ and a unitary flow $V$ on $E\Hil$ such that
$\|[E,x]\|\approx0$ for $x\in\G$, $\|(1-E)\xi\|\approx0$ for $\xi\in
\omega$ and $\|E\alpha_t(x)E-V_tExEV_t^*\|\leq \eps\|x\|$ for
$x\in \G$ and $t\in [-1,1]$.
Set $\beta_t=\Ad\,V_t$ on $\B(E\Hil)$.
From the above expression for $U^{(h,\alpha)}_t$ we note that
$$
Ev_t=EwU^{(h,\alpha)}\alpha_t(w^*)\approx
EwEU^{(EhE,\beta)}_t\beta_t(EwE)^*
$$
for $t\in [-1,1]$.
Thus  replacing $EwE$ by a close unitary in
$\B(E\Hil)$  we obtain a $\beta$-cocycle $b$ in $\B(E\Hil)$ such
that $Ev_t\approx Eb_t$ for $t\in [-1,1]$.
It follows that
$E\Ad\,v_t\alpha_t(x)E\approx \Ad\,b_t\beta_t(ExE)$ for $x\in\F$ and
$t\in [-1,1]$.
Since $\Ad\,b_t\beta_t(EeE)\approx
E\Ad\,b_t\alpha_t(e)E\approx E\Ad\,v_t(e)E\approx
E\Ad\,u_t(e)E\approx EeE$ for $t\in [-1,1]$, by Lemma \ref{pert1},
there is a subprojection $F$ of $E$ such that $Fp\approx Ep$,
$Fe\approx F$  and $\|[F,H']\|\approx 0$, where $H'$ is the
self-adjoint generator of $t\mapsto b_tV_t$.

Since $\|x-pxp\|\approx0$ for $x\in\F$, we have
$\|[F,x]\|\approx \|[E,x]\|\approx0$ for $x\in\F$ and since
$p\xi\approx \xi$ for $\xi\in\omega$ we have
$\|(1-F)\xi\|\approx \|(1-E)\xi\|\approx0$ for $\xi\in \omega$.
We conclude that $F\Ad\,u_t\alpha_t(x)F\approx
F\Ad\,v_t\alpha_t(x)F\approx F\Ad\,e^{itH'}(ExE)F\approx
\Ad\,e^{itFH'F}(FxF)$ for $x\in \F$.
\end{pf}

\begin{lem}
Let $A$ be a non-unital \cstar\ and $\alpha$ a flow on $A$.
Let $u$ be an $\alpha$-cocycle in $M(A)$.
Then there exists an approximate
identity $(e_\mu)_{\mu\in I}$ in $A$ such that
$\max\{\|\alpha_t(e_\mu)-e_\mu\|:\ t\in [-1,1]\}$ and
$\max\{\|[e_\mu,u_t]\|:\ t\in [-1,1]\}$ converge to zero.
Moreover
one may assume that there is another approximate identity
$(p_\mu)_{\mu\in I}$ with the same index set $I $ satisfying the
same conditions as $(e_\mu)$  and $e_\mu p_\mu=p_\mu$ for $\mu\in I$.
\end{lem}
\begin{pf}
Define a flow $\gamma$ on $M_2\otimes A$ by $\gamma_t(e_{12}\otimes
x)=e_{12}\otimes\alpha_t(x)u_t^*$ for $x\in A$. (Thus
$\gamma_t(e_{11}\otimes x)=e_{11}\otimes \alpha_t(x)$ and
$\gamma_t(e_{22}\otimes x)=e_{22}\otimes \Ad\,u_t\alpha_t(x)$.) We
choose an approximate identity $(f_\mu)$ in $M_2\otimes A$ such that
$\max\{\|\gamma_t(f_\mu)-f_\mu\|:\ t\in [-1,1]\} \to0$. By taking a
net in the convex combinations of $\{f_\mu\}$ we may further suppose
that $\|[e_{ij}\otimes 1,f_\mu]\| \to0$. Then we define $e_\mu\in A$
by
$$
1\otimes e_\mu={1\over 2}\sum_{i}(e_{i1}\otimes 1) f_\mu (e_{1i}\otimes
1),
$$
which is almost equal to $f_\mu$. Thus it follows that
$\|\gamma_t(1\otimes e_\mu)-1\otimes e_\mu\|\leq\|\gamma_t(1\otimes
e_\mu)-\gamma_t(f_\mu)\| +\|\gamma_t(f_\mu)-f_\mu\|+\|f_\mu-1\otimes
e_\mu\|$, which converges to zero uniformly in $t$ on $[-1,1]$.
Since $\gamma_t(1\otimes e_\mu)=e_{11}\otimes
\alpha_t(e_\mu)+e_{22}\otimes u_t\alpha_t(e_\mu)u_t^*$, this
completes the proof for the first part. To prove the additional
assertion, we choose two continuous functions $f,g$ from $[0,1]$
onto $[0,1]$ such that $f(0)=g(0)=0$, $f(1)=g(1)=1$ and $fg=g$.
Then the pair $f(e_\mu)$ and $g(e_\mu)$ satisfy
$f(e_\mu)g(e_\mu)=g(e_\mu)$.
One can prove that $f(e_\mu)$ (resp.\ $g(e_\mu)$) is an approximate
identity satisfying the required
properties.
\end{pf}

\begin{cor}\label{Her}
Let $\alpha$ be a flow on $A$ and let $B$ be an $\alpha$-invariant
hereditary \cstarsub\ of $A$ such that $B$ generates $A$ as an
ideal.
Then $\alpha$ is quasi-diagonal $($resp.\ pseudo-diagonal$)$ if
and only if $\alpha|_B$ is quasi-diagonal $($resp.\  pseudo-diagonal$)$.
\end{cor}
\begin{pf}
The ``only if'' part follows from the definition even if $B$ is an
arbitrary $\alpha$-invariant \cstarsub\ of $A$.
(See also Theorems~\ref{Ch} and \ref{Ch'}.)

Suppose that $\alpha|_B$ is quasi-diagonal (resp.\  pseudo-diagonal).
If $A$ is separable (or has a strictly positive element), then
$B\otimes \K$ and $A\otimes \K$ are isomorphic with each other (see
\cite{LB77}), where $\K$ is the separable \cstar\ of compact
operators. Under this identification, $\alpha\otimes\id$ on
$A\otimes\K$ is a cocycle perturbation of $\alpha|_{B\otimes \id}$ (see
\cite{K06}). Thus the ``if'' part follows from Propositions~\ref{pert}
and \ref{PS} in the separable case.

Suppose that $A$ is not separable.
Let $\F$ be a finite subset of $A$.
Since the linear span of $ABA$ is dense in $A$, there is a
countable subset $\G$ of $AB$ such that the closed linear span of
$\{xy^*:\ x,y\in\G\}$ contains $\F$.
Let $A_1$ be the $\alpha$-invariant \cstarsub\ of $A$ generated by $\G$.
Then $A_1\supset \F$.
Since $\alpha_s(x)^*\alpha_t(y)\in
A_1\cap B$ for $x,y\in \G$, the hereditary \cstarsub\ $B_1=A_1\cap
B$ of $A_1$ is essential, i.e.\  it generates $A_1$ as an ideal of
$A_1$.
Since $\alpha|_{B_1}$ is quasi-diagonal (resp.\ pseudo-diagonal),
it follows that $\alpha|_{A_1}$  is quasi-diagonal (resp.\
pseudo-diagonal).
Since $\F$ is arbitrary         this completes the proof.
\end{pf}

Recall that  pseudo-diagonality follows from quasi-diagonality.

\begin{prop}\label{KMS}
Suppose that $\alpha$ is a pseudo-diagonal flow on a unital \cstar\
$A$. Then $\alpha$ has a KMS state for all inverse temperatures
including $\pm\infty$.
\end{prop}
\begin{pf}
Let $\F$ be a finite subset $\F$ of $A$ and $\eps>0$. For each
$(\F,\eps)$ we have a flow $\beta$ on a finite-dimensional
\cstar\ $B$ and a CP map $\phi$ of $A$ into $B$ such that
$\phi(1)=1$, $\|\phi(x)\|\geq (1-\eps)\|x\|$ and
$\|\phi(x)\phi(y)-\phi(xy)\|\leq \eps\|x\|\|y\|$ for $x,y\in \F$
and $\|\beta_t\phi(x)-\phi\alpha_t(x)\|\leq \eps\|x\|$ for $x\in
\F$ and $t\in [-1,1]$.
Here we have replaced the condition
$\|\phi\|\leq1$ by $\phi(1)=1$ since $A$ is unital.
To justify this
we note that we may assume that $1\in\F$, which entails that
$\|\phi(1)^2-\phi(1)\|\leq \eps$ and
$\|\beta_t(\phi(1))-\phi(1)\|\leq \eps$ for $t\in [-1,1]$.
By functional calculus for small $\eps$  we obtain a projection $p$
from $\phi(1)$.
Since $\|\beta_t(p)-p\|$ is of order $\eps$ for
$t\in[-1,1]$  we can perturb $\beta$ by a $\beta$-cocycle which
differs from $1$ on $[-1,1]$ by up to order $\eps$ and
suppose that $\beta_t(p)=p$.
Replacing $B$ by $pBp$ and $\phi$ by
$q\phi(\,\cdot\,)q$ with $q=(p\phi(1)p)^{-1/2}$ and restricting
$\beta$  we can assume that $\phi$ is unital.
Since $\|q-p\|$ is of order $\eps$  we could start with a smaller
$\eps$ to obtain the right estimates.

There is a self-adjoint $h\in B$ such that $\beta_t=\Ad\,e^{ith}$.
Fix $\gamma\in \R$ and define a state $\varphi$ on $B$ by
$$
\varphi(Q)=\Tr(e^{-\gamma h}Q)/\Tr(e^{-\gamma h}),
$$
where $\Tr$ is a trace on $B$.
Then $\varphi$ is a KMS
state on $B$ with respect to $\beta$ at inverse temperature
$\gamma$.

Let $f_{(\F,\eps)}=\varphi\phi$ be a state on $A$  where
$\varphi$ and $\phi$ depend on $(\F,\eps)$.
Let $f$ be a
weak$^*$-limit point of $f_{(\F,\eps)}$, where the set $X$ of
$(\F,\eps)$ is a directed set in an obvious way. We fix a Banach
limit $\psi$ on $L^\infty(X)$ such that $f(x)$ is the $\psi$ limit
of $(\F,\eps)\mapsto f_{(\F,\eps)}(x)$ for $x\in A$.
Note that $f(x\alpha_t(y))$ is the $\psi$ limit of $(\F,\eps)\mapsto$
$\varphi(\phi(x\alpha_t(y)))$, which is close to $
\varphi(\phi(x)\beta_t\phi(y))$ around $\infty$.
Thus one can conclude that $f$ is a KMS state at $\gamma$.

A similar proof works for a KMS state for $\gamma=\pm\infty$ (or a
ground state and ceiling state).
See \cite{BR1,BR2} for more details
on KMS states.
\end{pf}

We may call such a state $f_{\F,\eps}$ on $A$ as above a {\em local
KMS state} (depending also on the choice of $B,\phi,\beta,h$ and
$\Tr$ on $B$) and a KMS state $f$ on $A$ obtained as a limit of
local KMS states {\em locally approximable}. It follows that the
locally approximable KMS states at an inverse temperature form  a
closed convex cone. It may be natural to ask whether all the KMS
states are locally approximable for a pseudo-diagonal flow on some
\cstar. An easy example of such will be given later.

We remind the reader that if $\alpha$ is approximately inner then we
obtain the same conclusion as in the above proposition \cite{BR1}.
The  proof is similar. Since there is a flow on a unital AF algebra
which has no KMS states for $\gamma>0$, we know that there is a
flow, on a unital AF algebra, which is not pseudo-diagonal. Obvious
examples of non-pseudo-diagonal flows are as follows:

\begin{example}
Let $\Omega$ be a compact Hausdorff space and  $\alpha$ a flow
of homeomorphisms of $\Omega$ such that no point of $\Omega$ is
fixed under $\alpha$.
We denote by the same symbol $\alpha$ the flow
of the \cstar\ on $A=C(\Omega)$ which naturally arises as
$\alpha_t(f)(\omega)=f(\alpha_{-t}(\omega))$ for $ f\in C(\Omega)$ and
$\omega\in \Omega$.
Then the flow $\alpha$ is not pseudo-diagonal since
if $\alpha$ has a KMS state for non-zero inverse
temperature then $\alpha$ acts trivially on $\pi(A)''$, where $\pi$
is the associated GNS representation of $A$, $($since $\pi(A)''$ is
commutative$)$ and this  implies the existence of fixed points under
$\alpha$ in $\Omega$.
\end{example}

\begin{example}\label{ex2}
Define a flow $\alpha$ on the \cstar\ $C_0(\R)$ by
$\alpha_t(f)(s)=f(s-t)$.
Then $\alpha$ is not pseudo-diagonal.
If one defines self-adjoint operators $P$ and $Q$ by
$P\xi(s)=-i{d\over {ds}}\xi(s)$ and $Q\xi(s)=s\xi(s)$ then  there is a finite sequence $f_1,f_2,\ldots,f_n$ in
$C_0(\R)$ and $\eps>0$ such that if a finite-rank projection $E$ on
$L^2(\R)$ satisfies  $\|Ef_i(Q)E\|\geq (1-\eps)\|f_i\|$ and
$\|[E,f_i(Q)]\|\leq \eps\|f_i\|$ for $i=1,\ldots,n$  then
$\|[E,P]\|>\eps$.
$($This statement appears considerably stronger than
the statement that if $E$ is a finite-rank projection on $L^2(\R)$ such
that $\|E\Omega_0\|>1/2$ then $\|[E,Q]\|+\|[E,P]\|>1/3$, where
$\Omega_0=\pi^{-1/4}e^{-s^2/2}$ is the {\em vacuum} vector.$)$
\end{example}

First consider the paranthetic assertion
and note that
$(P-iQ)(P+iQ)\geq 1$ and $(P-iQ)\Omega_0=0$.
Assuming there is such
a projection $E$ with $\|E\Omega_0\|>1/2$  let $T=E(P+iQ)E$ and
$\gamma=\|[E,P]\|+\|[E,Q]\|<1$.
Then $T\xi=[E,P+iQ]\xi+(P+iQ)\xi$
for $\xi\in EL^2(\R)$ and this  implies that $\|T\xi\|\geq
(1-\gamma)\|\xi\|$.
Since $T^*E\Omega_0=E(P-iQ)E\Omega_0=E[P-iQ,E]\Omega_0$ we deduce that
$\|T^*E\Omega_0\|\geq \gamma$. Since $\|T^{-1}\|=\|(T^*)^{-1}\|$
(as operators on the finite-dimensional subspace $EL^2(\R)$) it follows
that $\|E\Omega_0\|/\gamma\leq (1-\gamma)^{-1}$ or $\gamma\geq
\|E\Omega_0\|/(\|E\Omega_0\|+1)>1/3$.
(This assertion is related to the Heisenberg uncertainty principle.)

To establish the principal assertion of Example~\ref{ex2}  we may
add an identity  to $C_0(\R)$, i.e.\  we may consider $\alpha$ as
acting on the continuous functions on $\R^+=\R\cup\{\infty\}$. We
define a unitary $u\in C(\R^+)$ by $u(t)=1$ for $|t|\geq 1$ and
$u(t)=e^{i\pi(t+1)}$ for $t\in (-1,1)$. Note that
$\alpha_t(u)=e^{-ib_t}u$ for $t\geq 0$, where $b_t$ is a continuous
function on $\R$ with $\supp(b_t)=[-1,1+t]$ such that $b_t(s)=1+s$
for $s\in [-1,-1+t]$, $b_t(s)=t$ for $s\in (-1+t,1)$ and
$b_t(s)=1+t-s$ for $s\in [1,1+t]$.

We fix $t_0\in (0,1/2)$ and $\eps\in (0,1/6)$ and introduce
$f,g,h\in C_0(\R)$ as in Lemma~\ref{pq} below. In particular it
follows that $f\alpha_t(g)=\alpha_t(g)$ and $fb_t=b_t$ for $t\in
[0,t_0]$, $(u-1)g=u-1$, and $b_{t_0}h=t_0h$.
By applying Theorem~\ref{Ch'} to $u,f,g,h, b_t\in [0,t_0]$ etc.\  we obtain a unital CP
map $\phi$ of $C(\R^+)$ into $M_n$ for some $n$ and a flow $\beta$
on $M_n$ such that $\phi(u^*)\phi(u)\approx1$,
$\phi(g)(\phi(u)-1)\approx\phi(u)-1\approx (\phi(u)-1)\phi(g)$,
$\phi(\alpha_t(g))\approx \beta_t(\phi(g))$  and
$\phi(\alpha_t(u))\approx \phi(e^{-ib_t})\phi(u)\approx
e^{-i\phi(b_t)}\phi(u)\approx \beta_t(\phi(u))$ for $t\in [0,t_0]$
in addition to the conditions which ensure the conclusion of  Lemma~\ref{pq}.
We construct the spectral projections $F,G,H\in M_n$ out of
$\phi(f),\phi(g),\phi(h)$ corresponding to $[1-\delta,1]$ with a
small $\delta>0$ as in Lemma~\ref{pq}.
In particular  this ensures that
$G(\phi(u)-1)\approx \phi(u)-1$ and $F\beta_t(G)\approx\beta_t(G)$
and $F\phi(b_t)\approx\phi(b_t)$ for $t\in [0,t_0]$.
By slightly modifying $F$ we  can suppose that $GF=G$.
By the polar decomposition of
$G\phi(u)G+1-G\approx \phi(u)$ we obtain a unitary $W\in M_n$ such
that $W=GWG+1-G$ and  $e^{-i\phi(b_t)}W\approx\beta_t(W)$.
Let $V$ be a unitary flow in $M_n$ such that $\beta_t=\Ad\,V_t$.
Since $FV_tGF\approx V_tG$  there is a unitary
$Y_t\in FM_nF$ such that $Y_tG\approx V_tG$. We may suppose that
$t\in [0,t_0]\mapsto Y_t$ is continuous with $Y_0=F$.
Since $W(F-G)=
F-G$ and $F\beta_t(G)\approx \beta_t(G)$  we deduce  that
$F\beta_t(W)F\approx Y_tWY_t^*$ where  $W$ is now regarded as a
unitary in $FM_nF$.
By using $Fe^{-i\phi(b_t)}F\approx
e^{-iF\phi(b_t)F}$ we thus deduce that $Y_tWY_t^*W^*\approx
e^{-iF\phi(b_t)F} $ in $M_n$ for $t\in [0,t_0]$.
Hence there is a
self-adjoint $d_t\in FM_nF$ such that $d_t\approx 0$ and
$Y_tWY_t^*W^*=e^{-iF\phi(b_t)F}e^{id_t}$, where $t\mapsto d_t$ is
continuous.
Since $\det(Y_tWY_t^*W^*)=1$ we obtain
$-\Tr(F\phi(b_t)F)+\Tr(d_t)\in 2\pi\Z$.
Since $Y_tWY_t^*W^*=1$, $b_t=0$ and $d_t=0$ at $t=0$ it follows that
$\Tr(F\phi(b_t)F)=\Tr(d_t)$.
Note that $\Tr(F\phi(b_t)F)\geq
t\Tr(F\phi(h)F)\geq t(1-\delta)\Tr(FHF)$ and $\Tr(FHF)$ is almost
greater than $\dim(F)/3$ (by Lemma~\ref{pq}). Since
$|\Tr(d_t)|\leq\|d_t\|\dim(F)\approx0$  this gives a contradiction
for some $t$ away from $0$. (See \cite{EL,BEEK} for more details.)
{\hfill$\square$\vspace{3mm}}

%\medskip

\begin{lem}\label{prev}
For any $\eps_1,\eps_2,\eps'>0$ there is a $C>0$ such that the following holds.

For any $h,p\in A_{sa}$ such that $0\leq h\leq 1$, $0\leq p\leq 1$, and
$\|hp-p\|<\delta$  the spectral projections $F$ of $h$ and $G$ of
$p$ corresponding to $[1-\eps_1,1]$ and $[1-\eps_2,1]$,
respectively, satisfy $\|FG-G\|<\eps'+C\delta$.
\end{lem}
\begin{pf}
The proof is similar to that of Lemma \ref{proj}.
\end{pf}

\begin{lem}\label{pq}
Fix $t_0\in (0,1/2)$ and $\eps\in (0,1/6)$.
Define $g\in
C_0(\R)$ by $g(s)=0$ for $|t|>1+\eps$ and $g(s)=1$ for $t\in [-1,1]$
and by linearity elsewhere.
Define $f\in C_0(\R)$ by
$f=\alpha_{-\eps}(g)\vee \alpha_{t_0+\eps}(g)$ where $\alpha$ is
the translation flow of Example~$\ref{ex2}$.
$($Note that $f(s)=1$ for $s\in [-1-\eps,1+t_0+\eps]$ and
$f\alpha_t(g)=\alpha_t(g)$ for $t\in [0,t_0]$.$)$
Define $h\in C_0(\R)$ by $h(t)=0$ for $t<-1+t_0$
and $t>1$, $h(t)=1$ for $t\in [-1+t_0+\eps,1-\eps]$ and by linearity
elsewhere.
$($Note that $f\cdot\bigvee\{\alpha_s(h):\ |s|\leq t_0+3\eps\}=f$.$)$

Let $\delta\in (0,1)$ and suppose that $\alpha$ is pseudo-diagonal.
Then for any $\eps'>0$ there is a unital CP map $\phi$ of
$C_0(\R^+)$ into $M_n$ and a flow $\beta$ on $M_n$ satisfying
the following assertion.
If $F$, $G$, and $H$ are spectral projections of $\phi(f)$, $\phi(g)$,
and $\phi(h)$, respectively, corresponding to $[1-\delta,1]$, then
$\|F\beta_t(G)-\beta_t(G)\|<\eps'$ for $t\in [0,t_0]$ and $\dim
F\leq 3\dim H$.
\end{lem}
\begin{pf}
The estimate  $\|FG-G\|<\eps'$ follows from Lemma~\ref{prev}  by assuming
$\|\phi(f)\phi(g)-\phi(g)\|\approx0$.
If $\phi(\alpha_t(g))\approx\beta_t(\phi(g))$ sufficiently closely then
$\beta_t(G)$ is almost dominated by the spectral projection
$G_t$ of $\phi\alpha_t(g)$ corresponding to $[1-3\delta/2,1]$ (and
almost dominates the spectral projection corresponding to
$[1-\delta/2,1]$). (See Lemma 2.2 of \cite{BEEK}.)
If $\phi(f)\phi\alpha_t(g)\approx \phi\alpha_t(g)$ sufficiently closely
then  $FG_t\approx G_t$.
Thus, assuming
$\|\phi(f)\phi\alpha_t(g)-\phi\alpha_t(g)\|\approx0$, it follows
that $F\beta_t(G)\approx\beta_t(G)$ for $t\in [0,t_0]$.

Let $t_1=-t_0-3\eps$, $t_2=2-2t_0-6\eps$, and $t_3=4-3t_0-9\eps$ and
note that $\alpha_{t_1}(h)\vee\alpha_{t_2}(h)\vee
\alpha_{t_3}(h)\cdot f=f$ and that there are non-negative
$f_1,f_2,f_3\in C_0(\R)$ such that $f=f_1+f_2+f_3$ and
$\alpha_{t_i}(h)f_i=f_i$.
Suppose that $\dim F>3\dim H$ however we
choose $\phi$ and $\beta$.
Then, since $\dim F>\dim(\beta_{t_1}(H)\vee\beta_{t_2}(H)\vee \beta_{t_3}(H))$,
there is a state $\varphi$ on $M_n$ such that $\varphi(F)=1$ and
$\varphi(\beta_{t_i}(H))=0$ for $i=1,2,3$.
By assuming that $\|\phi(f_i)\phi\alpha_{t_i}(h)-\phi(f_i)\|\approx0$ etc.\  we would
have  $\varphi(\phi(f_i))\approx0$ for $i=1,2,3$, which implies
that $\varphi(\phi(f))\approx 0$.
But since $\varphi(\phi(f))\geq 1-\delta$ due to $\varphi(F)=1$
this is a contradiction.
Hence  $\dim F\leq 3\dim H$ follows if $\phi(f_i)\phi\alpha_{t_i}(h)\approx
\phi(f_i)$.
\end{pf}

\begin{example}
Let $\Dk$ denote the unit disk $\{z\in\C:\ |z|\leq1\}$. Define a
flow $\alpha$ of homeomorphisms of $\Dk$ by $\alpha_t(z)=ze^{it}$.
Then  the induced  flow on $C(\Dk)$ is quasi-diagonal. More
generally let $v$ be a continuous function on $[0,1]$ of finite
variation and define a flow $\alpha'$ on $\Dk$ by
$\alpha'_t(z)=e^{itv(|z|)}z$. Then the induced flow on $C(\Dk)$ is
quasi-diagonal. Note that the origin is a fixed point which is
neither absorbing nor repelling.
\end{example}

We shall prove the first assertion here. The second one will not be
proved but follows from the proof of Proposition \ref{hom2} given
later.

Let $\alpha$  denote the induced flow on  $C(\Dk)$
and  $\beta$  the rotation flow on $C(\T)$, i.e.\
$\beta_t(x)(z)=x(ze^{-it})$, where $\T=\{z\in\C:\ |z|\leq1\}$.
We regard $C(\T)$ as acting on $L^2(\T)$.
Then $\beta$ is implemented
by the unitary flow $U$ defined by $U_t\xi(z)=\xi(ze^{-it})$.
Note that $U_t=\sum_{k=-\infty}^\infty e^{ikt}P_k$ where $P_k$ is a
rank-one projection.
For $r\in [0,1]$ let $\pi_r$ be the restriction
map of $C(\Dk)$ onto $C(\T)\colon \pi_r(x)(z)=x(rz)$.
For $n\in\N$, we define a covariant representation $\rho_n$ of $C(\Dk)$ by
$\rho_n=\bigoplus_{k=0}^n \pi_{k/n}$ with the unitary flow $U^{(n)}$
defined by $U^{(n)}_t=\bigoplus_{k=0}^n U_{t}$.

Let $\F$ be a finite subset of $C(\Dk)$ and $\eps>0$.
For any $\delta>0$ there is an $n\in \N$ such that $\|\rho_n(x)\|\geq
(1-\delta)\|x\|$ for all $x\in \F$ and
$\|\pi_r(x)-\pi_s(x)\|\leq\delta\|x\|$ if $x\in\F$ and $\|r-s\|\leq 1/n$.
We find a decreasing sequence
$T_0=F_0,G_0,T_1,F_1,G_1,\ldots,T_n,F_n,G_n=0$ of non-negative
operators in the convex hull of the $P_k$ such that all $F_k$ and
$G_k$ are projections, $\|(F_k-G_k)\pi_{k/n}(x)(F_k-G_{k})\|\geq
(1-\delta)\|\pi_{k/n}(x)\|$ for $x\in\F$ and
$\|[T_k,\pi_{k/n}(x)]\|\leq\delta\|x\|$ for $x\in\F$.
We construct the sequence in the  reverse order.

After choosing $G_k$, since $\pi_{k/r}(x)$ is not compact, one can
choose $F_k$ to satisfy the condition
$\|(F_k-G_k)\pi_{k/n}(x)(F_k-G_k)\|\geq (1-\delta)\|\pi_{k/n}(x)\|$.
By the general theory of quasi-central approximate units, (see
\cite{Ar} or \cite{Ped}), one can choose $T_k\geq F_k$. If $T_{k}$
is chosen we set $G_{k-1}$ to be the support projection of $T_{k}$.
After repeating this process  a finite number of times we construct
$F_0$. Since $\|[T_k,\pi_{k/n}(x)]\|\leq\delta\|x\|$ is void for
$k=0$ we can set $T_0=F_0$.

We define a finite-rank projection $E=S^*S$ on $\bigoplus_{k=0}^nL^2(\T)$
with
$S=((T_0-T_1)^{1/2},(T_1-T_2)^{1/2},\ldots,(T_n-T_{n+1})^{1/2})$
with $T_{n+1}=0$.
Since $SS^*=T_0$, $E$ is indeed a finite-rank
projection.
Since all the  $T_k$ commute with $U$, it follows that
$[E,U^{(n)}_t]=0$.
Since $E$ is tri-diagonal, $[E,\rho_n(x)]$ is
expressed as the sum of the diagonal part $\bigoplus_{k=0}^n
[T_k-T_{k+1},\pi_{k/n}(x)]$, the upper off-diagonal part
$$
\bigoplus_{k=0}^{n-1}\{E_{k,k+1}\pi_{(k+1)/n}(x)
-\pi_{k/n}(x)E_{k,k+1}\},
$$
and the lower off-diagonal part
$$
\bigoplus_{k=0}^{n-1}\{E_{k+1,k}\pi_{k/n}(x)
-\pi_{(k+1)/n}(x)E_{k+1,k}\},
$$
where
$E_{k,k+1}=E_{k+1,k}=(T_k-T_{k+1})^{1/2}(T_{k+1}-T_{k+2})^{1/2}
=(T_{k+1}-T_{k+1}^2)^{1/2}$.
Thus one can conclude that $\|[E,\rho_n(x)]\|\leq\eps\|x\|$ for
$x\in \F$ for a sufficiently small $\delta$. The diagonal part of
$E\rho_n(x)E$ is given by the direct sum of
$$
(T_k-T_{k+1})\pi_{k/n}(x)(T_k-T_{k+1})
+E_{k,k+1}\pi_{(k+1)/n}(x)E_{k+1,k}+E_{k,k-1}\pi_{(k-1)/n}(x)E_{k-1,k},
$$
over $k=0,1,\ldots,n$, where $E_{k,k+1}$ etc.\  are given above and
the term should be omitted if $k+1=n+1$ or $k-1=-1$.
Hence the norm of this is greater than or equal to
$$
\bigoplus_{k=0}^n (F_k-G_k)\pi_{k/n}(x)(F_k-G_k).
$$
Thus we obtain  $\|E\rho_n(x)E\|\geq(1-\eps)\|x\|$ for $x\in\F$
for a small $\delta$.
{\hfill$\square$\vspace{3mm}}

Since $\pi_0(x)=x(0)$, $\pi_0$ is an $\alpha$-invariant character.
Hence we could choose $t\mapsto 1$ for a unitary flow implementing
$\alpha$ instead of the $U$ which has spectrum $2\pi\Z$, but then
the above proof would fail.

The above proof was taken from the proof of Proposition~3 of \cite{Voi91}.
It is appropriate to indicate how to prove Theorem
\ref{hom} at this point.
First we establish the following analogue of
Proposition~3 of \cite{Voi91}.

\begin{prop}\label{hom1}
Let $\alpha$  $($resp.\  $\beta$$)$ be a flow on a \cstar\ $A$ $($resp.\  $B$$)$.
Let $\{\phi_s;\ s\in [0,1]\}$ be a homotopy of homomorphisms of $A$
into $B$ such that $\phi_s\alpha_t=\beta_t\phi_s$ and
$\bigcap_s\Ker(\phi_s)=\{0\}$.
If $\beta|_{\phi_1(A)}$ is
quasi-diagonal  then $\alpha$ is quasi-diagonal.
\end{prop}
\begin{pf}
Let $(\rho,V)$ be a covariant representation of $(B,\beta)$ such
that $\rho\times V$ is faithful and contains no non-zero compact
operators in its range. Then, by the assumption,
$\bigoplus_s\rho\phi_s$ is faithful and $(\rho\phi_1(A),V)$ is
quasi-diagonal. Let $H$ denote the self-adjoint generator of $V$.
Let $\F$ be a finite subset of $\A$ and $\eps>0$. There is a
self-adjoint compact operator $K$ on $\Hil_\rho$ such that
$\|K\|<\eps/2$ and $H_1=H+K$ is diagonal. For any small constant
$\delta>0$, there is an $n\in \N$ such that if $|s_1-s_2|\leq 1/n$
and $x\in \F$ then $\|\phi_{s_1}(x)-\phi_{s_2}(x)\|\leq \delta\|x\|$
and if $x\in\F$ then $\max_k\|\pi_{k/n}(x)\|\geq(1-\delta)\|x\|$.
There is a finite increasing sequence $G_0=0,
F_0,T_0,G_1,F_1,T_1,G_2,\ldots,G_n,F_n=T_n$ of non-negative compact
operators in the maximal commutative von Neumann algebra generated
by a family of minimal projections commuting with $H_1$ such that
all $G_k,F_k$ are projections,
$\|(F_k-G_k)\pi_{k/n}(x)(F_k-G_k)\|\geq (1-\delta)\|\pi_{k/n}(x)\|$
for $x\in\F$  and $\|[T_k,\phi_{k/n}(x)]\|<\delta\|x\|$ for
$x\in\F$. Let $\pi=\bigoplus_{k=0}^n\rho\phi_{k/n}$ and
$U=\bigoplus_{n=0}^n V$. In this covariant representation space
$\bigoplus_{k=0}^n \Hil_\rho$ we define a finite-rank projection $E$
as $S^*S\in M_{n+1}\otimes\B(\Hil_\rho)$, where $S$ is the row
vector $(T_0, (T_1-T_0)^{1/2},\ldots, (T_n-T_{n-1})^{1/2})$. If
$\delta$ is sufficiently small, one can show that $(\pi,U)$ and $E$
satisfy Condition (3) of Theorem~\ref{Ch} for $(\F,\eps)$. (See the
proof of Proposition~3 of \cite{Voi91} for more details.)
\end{pf}

Theorem \ref{hom} follows from  Proposition~\ref{hom1} exactly as in
Theorem 5 of \cite{Voi91}.
Let us reproduce the proof here.
We have two flows $\alpha$ on $A$ and $\beta$ on $B$ such that
$(B,\beta)$ dominates $(A,\alpha)$, i.e.\  there are intertwining
homomorphisms $\phi\colon A \to B$ and $\psi\colon B \to A$ such that $\psi\phi$
is homotopic to the identity in the endomorphisms of $A$ commuting with
$\alpha$.
The assumption that $\beta$ is quasi-diagonal implies that
$\beta|_{\phi(A)}$ is quasi-diagonal and hence that the flow
$\dot{\alpha}$ on $A/\Ker(\phi)$ induced from $\alpha$ is
quasi-diagonal.
Let $D=A\oplus A/\Ker(\phi)$ with the flow
$\alpha\oplus\dot{\alpha}$ and consider two intertwining
homomorphisms $\psi\phi\oplus \pi$ and $\id\oplus \pi$ from $A$ into
$D$, where $\pi$ is the quotient map of $A$ onto $A/\Ker(\phi)$.
Then, since $\Ran(\psi\phi\oplus\pi)$ is isomorphic to
$A/\Ker(\phi)$, we conclude that
$\alpha\oplus\dot{\alpha}|_{\Ran(\psi\phi\oplus\pi)}$ is
quasi-diagonal.
Since $\psi\phi\oplus\pi$ is homotopic to
$\id\oplus\pi$ in the intertwining homomorphisms and $\id\oplus\pi$
is injective Proposition~\ref{hom1} implies that $\alpha$ is
quasi-diagonal.
This concludes the proof of \ref{hom}.
(See \cite{Voi91} for another formulation.)

We can also show the following variant of Proposition~3 of
\cite{Voi91}.

\begin{prop}\label{hom2}
Let $\alpha$ be a flow on a quasi-diagonal \cstar\ $A$  and define a
flow $\beta$ on $B=A\otimes C[0,1]$ by
$\beta_t(x)(s)=\alpha_{st}(x(s))$. Then $\beta$ is quasi-diagonal.
\end{prop}
\begin{pf}
Let $(\pi,U)$ be a covariant representation of $(A,\alpha)$ such
that $\pi\times U$ is faithful. For $s\in [0,1]$ we define a map
$\phi_s$ of $B$ onto $A$ by $\phi_s(x)=x(s)$. For each $n\in \N$ let
$L_n=1+1/2+1/3+\ldots +1/n$ and let $s_k=L_n^{-1}\sum_{m=1}^k1/m$
for $k=1,2,\ldots,n$  with $s_0=0$. We define a representation
$\pi_n$ of $B$ by $\pi_n=\bigoplus_{k=0}^n\pi\phi_{1-s_k}$ and a
unitary flow  $U^{(n)}_t=\bigoplus_{k=0}^n U_{(1-s_k) t}$ which
implements $\beta$. Given $\F$ and $\eps>0$ we construct the
required finite-rank projection in the representation space
$\Hil_n=\bigoplus_{k=0}^n\Hil_\pi$ of $\pi_n$ for some large $n$.
The finite-rank projection $E=(E_{k\el})$ is defined just as before
as the tri-diagonal matrix $S^*S$ by choosing a finite increasing
sequence $T_0,T_1,\ldots,T_n$ of finite-rank non-negative operators
as above; $E_{kk}=T_k-T_{k-1}$ and
$E_{k,k+1}=E_{k+1,k}=(T_k-T_{k-1})^{1/2}(T_{k+1}-T_{k})^{1/2}=(T_k-T_k^2)^{1/2}$.
Since $U^{(n)}$ is not the direct sum of $n+1$ copies of the same
flow  it is not sufficient to assume $T_k$  almost commutes  with
$U_{t(1-s_k)}$. To achieve  $\Ad\,U_t^{(n)}(E)\approx E$  we must
also have
$$
U_{(1-s_k)t}E_{k,k+1}U_{(1-s_{k+1})t}^*\approx E_{k,k+1}.
$$
Since $E_{k,k+1}$ almost commutes with $U$, this amounts to
$E_{k,k+1}U_{t(k+1)^{-1}L_n^{-1}}\approx E_{k,k+1}$, i.e.\
$(U_t-1)(T_k-T_k^2)\approx 0$ if $|t|\leq (k+1)^{-1}L_n^{-1}$.

Let $X$ be the linear subspace of $C_0(\R^+)$ consisting of
non-increasing $C^\infty$-functions $f$ on $\R$ such that $f(t)=1$
for all small $t$ and $f(t)=0$ for all large $t$.
We regard $H$ as
the generator of $t\mapsto \lambda_t$ in the multiplier algebra
$M(A\times_\alpha\R)$ and write $f(H)=\int \hat{f}(t)\lambda_tdt$.
Since $f(H)$ with $ f\in X$  is a subspace in $M(A\times_\alpha\R)$ and
contains the identity in its closure with respect to the strict
topology, there is an $f_-\in X$ such that
$\|[f_-(H),x]\|\leq(\eps/4)\|x\|$ for all $x\in \F$, where $\F$ is
the given finite subset of $A$.
Since this inequality is left
invariant under the dual flow  we may assume that $f_-(t)=1$ for
$t\leq 0$ and $f_-(t)<1$ for $t>0$. Let $L=\min\{t:\ f_-(t)=0\}$ and
define $f_+$ by $f_+(t)=1-f_-(t+L)$. On the other hand there is an
$L'>0$ such that
$$\|\chi_{(-L',L')}(H)x\chi_{(-L',L')}(H)\|\geq(1-\eps/2)\|x\|
$$
for $x\in\F$  where $\chi_{(-L',L')}$ is the characteristic function
of the interval $(-L',L')$ (and hence $\chi_{(-L',L')}(H)$ is an
open projection in the second dual of $A\times_\alpha\R$).
We define
$f_k\in C_0(\R)$ by
$$f_k(t)=f_+(t+kL+(2k+1)L')f_-(t-kL-(2k+1)L')
$$
for $k=0,1,2,\ldots$. We note that
$\supp(f_k)=[-(k+1)L-(2k+1)L',(k+1)L+(2k+1)L']$, $f_k(t)=1$ for
$t\in [-kL-(2k+1)L',kL+(2k+1)L']$, $f_kf_{k+1}=f_k$, and
$\|[f_k(H),x]\|\leq (\eps/2)\|x\|$ for $x\in\F$.

Let $P_0=\chi_{(-L',L')}(H)$  where $H$ now denotes the generator of $U$.
We choose a finite-rank non-negative operator $T_0$ on
$\Hil_\pi$ such that $T_0\leq f_0(H)$, $\|[T_0,U_t]\|<\eps$ for
$t\in [-1,1]$, and $\|[T_0,\pi(x)]\|\leq(\eps/2)\|x\|$ and
$\|P_0T_0\pi(x)T_0P_0\|\geq(1-\eps)\|x\|$ for $x\in \F$. This is
possible since the set of finite-rank operators $T$ satisfying
$0\leq T\leq f_0(H)$ forms a convex set invariant under $\Ad\,U_t$
and contains $f_0(H)$ in its closure with the strict topology on
$\K(\Hil_\pi)$ and $f_0(H)$ satisfies all the conditions required
for $T_0$ with stricter coefficients. Let $G_1$ be the support
projection of $T_0$ and let
$P_1=\chi_{(-L-3L',-L-L')}(H)+\chi_{(L+L',L+3L')}(H)$. We will find
a finite-rank non-negative operator $T_1$ on $\Hil_\pi$ such that
$G_1\leq T_1\leq f_1(H)$, $\|[T_1,U_t]\|<\eps$ for $t\in [-1,1]$,
and $\|[T_1,\pi(x)]\|\leq(\eps/2)\|x\|$ and
$\|P_1T_1\pi(x)T_1P_1\|\geq(1-\eps)\|x\|$ for $x\in \F$. For this
purpose let $N$ be a large number and let $F$ be a finite-rank
projection such that $F\leq \chi_{(-L-L',L+L')}(H)$ and
$\|FU_tG_1U_t^*-U_tG_1U_t^*\|=\|U_t^*FU_tG_1-G_1\|\approx0$ for all
$t\in [-N,N]$. We then choose a finite-rank $T$ such that $F\leq
T\leq f_1(H)$, $\|[T,P_1]\|<\eps$, $\|[T,\pi\alpha_t(x)]\|\leq
(2\eps/3)\|x\|$ and $\|P_1T\pi\alpha_t(x)TP_1\|\geq
(1-2\eps/3)\|x\|$ for $x\in \F$ and $t\in [-N,N]$. Then by taking
$\int h(t)U_tTU_t^*dt$ instead of $T$ with an appropriate function
$h\geq0$ and assuming that $N$ is large enough we can see that all
the conditions are satisfied with $TG_1\approx G_1$ instead of
$G_1\leq T$. We then modify $T$ slightly to  obtain a $T_1$ which
satisfies all the conditions. Note that $T_1T_0=T_0$ and
$\|[T_1-T_0,\pi(x)]\|\leq \eps\|x\|$.

By repeating this process we obtain $T_0,T_1,\ldots,T_{n-1}$ such
that $T_k\leq f_k(H)$, $T_kT_{k-1}=T_{k-1}$, $\|[T_k,U_t]\|<\eps$
for $t\in [-1,1]$, $\|[T_k-T_{k-1},\pi(x)]\|\leq \eps\|x\|$ and
$\|P_kT_k\pi(x)T_kP_k\|\geq(1-\eps)\|x\|$ for $x\in\F$, where
$$P_k=\chi_{(-kL-(2k+1)L',-kL-(2k-1)L')}(H)
+\chi_{(kL+(2k-1)L',kL+(2k+1)L')}(H).
$$
Let
$P'_n=\chi_{(-\infty,-nL-(2n+1)L')}(H)+\chi_{(2nL+(2n+1)L',\infty)}(H)$.
By using that $A$ is quasi-diagonal, we then choose a finite-rank
projection $T_n$ such that $T_n\geq T_{n-1}$ and
$\|P'_nT_n\pi(x)T_nP'_n\|\geq (1-\eps)\|x\|$ and
$\|[T_n,\pi(x)\|\leq \eps\|x\|$ for $x\in\F$. If $|t|\leq
(k+1)^{-1}L_n^{-1}$, then $\|(U_t-1)(T_k-T_k^2)\|\leq
|L(k+1)+L'(2k+1)|(k+1)^{-1}L_n^{-1}\leq (L+2L')/L_n$. Since $L_n
\to\infty$ as $n \to\infty$ we obtain the desired sequence
$T_0,T_1,\ldots,T_n$ for some $n$.
\end{pf}

\begin{example}
Let $A_\theta$ denote the irrational rotation \cstar\ generated by
two unitaries $u,v$ satisfying $uv=e^{i2\pi\theta}vu$ with
$\theta\in (0,1)$ irrational.
$($Then $A_\theta$ is a unital simple AT-algebra with a unique tracial state.$)$
Let $\alpha$ be a flow on
$A_\theta$ such that $\alpha_t(u)=e^{it}u$ and $\alpha_t(v)=e^{i p
t}v$ with some $p\in\R$.
Then $\alpha$ is not pseudo-diagonal.
This follows because if $\omega$ is a KMS state for the inverse temperature
$\beta\not=0$ then one  must have
$e^{-\beta}=\omega(u^*\alpha_{i\beta}(u))=\omega(uu^*)=1$ which is
a contradiction.
Thus $\alpha$ has no KMS states and $\alpha$ is not
pseudo-diagonal.
\end{example}

\begin{prop}\label{AI}
Let $A$ be a quasi-diagonal \cstar\ and $\alpha$ an approximately
inner flow on $A$.
Then $\alpha$ is pseudo-diagonal.
\end{prop}
\begin{pf}
Suppose $A$ acts non-degenerately on a Hilbert space $\Hil$
such that $A$ is a quasi-diagonal set of $\B(\Hil)$.

Let $\F$ be a finite subset of $A$ and $\eps>0$.
By assumption there is an $h=h^*\in A$ such that
$\|\alpha_t(x)-\Ad\,e^{ith}(x)\|\leq \eps/3\|x\|$ for $x\in\F$
and $t\in[-1,1]$.
There is a finite-rank projection $E$ on $\Hil$
such that $\|ExE\|\geq (1-\eps)\|x\|$ and
$\|[E,x]\|\leq\eps\|x\|$ for $x\in\F$, and
$\|[E,h]\|<\eps/3$. Since $\|Ee^{ith}E-e^{itEhE}E\|<\eps/3$
for $t\in [-1,1]$, it follows that
$$
\|E\alpha_t(x)E-\Ad\,e^{itEhE}(ExE)\|\leq \eps\|x\|
$$
for all $x\in \F$. Note also that $\|ExEyE-ExyE\|\leq \eps
\|x\|\|y\|$ for $x,y\in \F$. By setting $B=\B(E\Hil)$,
$\beta_t=\Ad\,e^{itEhE}$ and $\phi(x)=ExE$ we obtain Condition (2)
of Theorem~\ref{Ch'}.~\end{pf}

If $\alpha$ is a pseudo-diagonal flow on $A$ and $B$ is an
$\alpha$-invariant \cstarsub\ of $A$ then $\alpha|_B$ is
pseudo-diagonal, i.e.\  pseudo-diagonality is preserved under passing
to an invariant \cstarsub.
But it is not evident that  this
property holds for approximate innerness.

\begin{prop}\label{AF}
If $\alpha$ is an AF flow  then $\alpha$ is quasi-diagonal.
\end{prop}
\begin{pf}
By  assumption the \cstar\ $A$ is an AF algebra and has an
increasing sequence $(A_n)$ of finite-dimensional C$^*$-subalgebras
such that $\bigcup_n A_n$ is dense in $A$ and $\alpha_t(A_n)=A_n$.
We choose a maximal abelian \cstarsub\ $D_n$ of $A_{n}\cap A_{n-1}'$
(with $A_0=\C1$) such that $\alpha$ is trivial on $D_n$ and let $D$
be the \cstarsub\ of $A$ generated by all $D_n$.
Let $(\phi_n)$ be a dense sequence in the characters of $D$.
Each $\phi_n$ uniquely extends to a pure $\alpha$-invariant state of $A$
which we also denote by $\phi_n$.
Note that $\bigoplus_n\pi_{\phi_n}$ is a faithful representation of~$A$.

In the GNS representation $\pi_\phi$ of $A$ for $\phi=\phi_n$ we
define a unitary flow $U$ by
$U_t\pi_\phi(x)\Omega_\phi=\pi_\phi\alpha_t(x)\Omega_\phi,\ x\in A$.
It follows that $(\pi_\phi,U)$ is a covariant representation of
$(A,\alpha)$. Denote by $E_n$ the finite-rank projection onto the
subspace $\pi_\phi(A_n)\Omega_n$.
Then $[E_n,\pi_\phi(x)]=0$ for $x\in
A_n$, $[E,U_t]=0$ and $\lim_n E_n=1$.
This shows that $(\pi_\phi(A),U)$ is quasi-diagonal.
Denoting $U$ by $U^\phi$ we conclude that
$(\bigoplus_n\pi_{\phi_n},\bigoplus_nU^{\phi_n})$ is
quasi-diagonal.
\end{pf}

From the above proof one can construct a projection of norm one
$\phi_n$ of $A$ onto $A_n$ such that
$\alpha_t|_{A_n}\circ\phi_n=\phi_n\circ\alpha_t$, from which follows
a stronger form of Condition (2) of Theorem \ref{Ch}. By using this
fact one can show that all the KMS states are locally approximable
for an AF flow (see 4.6.1 of \cite{Sak} and the comment after
Proposition \ref{KMS}). This remark also applies to approximate AF
flows \cite{K05}.

A quasi-diagonal flow on an AF algebra is not expected to be simply
a cocycle perturbation of an AF flow because there is a more general
type of AF flow (see \cite{Sak} for commutative derivations which
generate such flows). Specifically, there is a flow $\alpha$ on a
unital simple AF algebra $A$,  which is not a perturbation of an AF
flow, such that $A$ has an increasing sequence $(B_n)$ of
$\alpha$-invariant C$^*$-subalgebras of $A$ with dense union
satisfying  $\alpha|_{B_n}$ is uniformly continuous and $B_n\cong
A_n\otimes C[0,1]$ with $A_n$ finite-dimensional \cite{K01}. We take
an $\alpha$-invariant pure state $f$ such that $f|_{B_n}$ reduces to
an evaluation on $C[0,1]$. Then, in the GNS representation
associated with $f$, the subspace $\pi_f(B_n)\Omega_f$ is
finite-dimensional, which  gives the desired finite-rank
projections.

\begin{example}
Let $(A_n)$ be an increasing sequence of finite-dimensional
C$^*$-algebras and $A$ the closure of the union
$\bigcup_nA_n$. Let $B_n=A_1\oplus A_2\oplus\cdots\oplus A_n$. We
define an embedding of $B_n$ into $B_{n+1}$ as follows: if $k<n$
then $A_k$ of $B_n$ is identified with $A_k$ of $B_{n+1}$ and $A_n$
is mapped into $A_n\oplus A_{n+1}$ by duplication.
Let $B$ be the closure of the union $\bigcup_nB_n$ $($which is also defined as the
\cstar\ of bounded sequences $(x_n)$ with $x_n\in A_n$ such that
$\lim_nx_n$ converges in $A$$)$. We note that $B$ has many
finite-dimensional quotients. Let $I_n$ be the ideal generated by
all $A_i\subset B_m$ with $i\not=n\leq m$. Then $B/I_n$ is
isomorphic to $A_n$. Since $\bigcap I_n=\{0\}$, one concludes that
$B$ is quasi-diagonal. This shows that any flow $\beta$ on $B$ is
quasi-diagonal $($since $\beta$ fixes $I_n$$)$.

Let $\alpha$ be an approximately inner flow on $A$ and choose an
$h_n=h_n^*\in A_n$ such that $\Ad\,e^{ith_n}(x) \to \alpha_t(x)$ for
$x\in A$.
We can define a flow $\beta$ on $B$ as follows.
Let $\beta^{(n)}$ denote the flow on $B_n$ implemented by
$\bigoplus_{k\leq n}h_k$.
One shows that $\beta^{(n)}_t(x)$
converges for $x\in B$ and defines $\beta$ as the limit.

Conversely, if a flow $\beta$ is given on $B$  then we have an
$h_n=h_n^*\in A_n$ such that the induced flow on $B/I_n$ is given by
$Ad\,e^{ith_n}$ and can argue that $\Ad\,e^{ith_n}$ converges to a flow
on $A$.
Hence one concludes that $\beta$ is defined just as in the
previous paragraph.
\end{example}

\begin{example}
Let $A$ be a residually finite-dimensional \cstar\ and  $\alpha$
a flow on $A$ which fixes each ideal of $A$.
Then $\alpha$ is quasi-diagonal.
This follows because $A$ has a separating family of
finite-dimensional representations which must be covariant under
$\alpha$.
Thus the direct sum of these representations gives the
required faithful representation of~$A$.
\end{example}

\begin{example}\label{GICAR}
Let $\gamma$ denote the periodic flow on the UHF algebra
$\bigotimes_{n=1}^\infty M_{2}$ of type $2^\infty$ given by
$\gamma_t=\bigotimes_n\Ad(1\oplus e^{2\pi i t})$ and let $A$ be the
fixed point algebra of $\gamma$. The dimension group of $A$ is
isomorphic to $\Z[t]$ with the positive cone of strictly positive
functions on the open interval $(0,1)$. There is a decreasing
sequence $I_1,I_2,\ldots$ of ideals of $A$ such that
$A/I_1\cong\C1$, $\bigcap_n I_n=\{0\}$, and $I_n/I_{n+1}\cong \K$
for $n\geq 1$, where $\K$ is the compact operators on a separable
infinite-dimensional Hilbert space. It follows from the next lemma
that any flow on $A/I_n$ is quasi-diagonal and this  then implies
that any flow on $A$ is quasi-diagonal.
\end{example}

\begin{lem}
If $A$ is a type I AF algebra then any flow on $A$ is quasi-diagonal.
\end{lem}
\begin{pf}
There is a strictly increasing family $\{I_\mu\}$ of (closed) ideals
of $A$ indexed by $\mu$ in a segment $[0,\nu]$ of ordinals such that
$I_0=\{0\}$, $I_\nu=A$, $\bigcup_{\mu<\gamma}I_\mu$ is dense in
$I_\gamma$ for any limit ordinal $\gamma$, and $I_{\mu+1}/I_\mu$ is
generated as an ideal by a minimal projection for any $\mu<\nu$,
i.e.\  $I_{\mu+1}/I_\mu\cong \K$ or otherwise $M_m$ for some
$m=1,2,\ldots$ (see, e.g.\  \cite{Ped} for type I \cstars).
In the following we allow  $I_{\mu+1}/I_\mu\cong 0$ and  call this a
composite series for $A$.
Note that any flow $\alpha$ on $A$ fixes
each ideal (because the ideal is generated by projections). We shall
prove the statement that any flow on $A$ (with a composite series
indexed by $[0,\nu]$) is quasi-diagonal by induction on $\nu$.
If $\nu=1$ this is obvious since $A\cong\K$ or $M_m$ for some $m$ and
any flow on $A$ is inner.
(If $A\cong \K$ we use the Weyl-von Neumann theorem.)

Suppose that this is shown for any $\nu<\sigma$. Let $A$ be a type I
AF algebra with a composite series $\{I_\mu\}$ with $\mu\in
[0,\sigma]$.
If $\sigma$ is a limit ordinal, then
$\bigcup_{\nu<\sigma}I_\nu$ is dense in $A$. Since $I_\nu$ has a
composite series indexed by $[0,\nu]$  the induction hypothesis
implies that any flow on $I_\nu$ is quasi-diagonal  which in turn
implies that any flow on $A$ is quasi-diagonal.
If $\sigma$ is not a
limit ordinal (and $I_{\sigma-1}\not=A$) then there is a minimal
projection in $A/I_{\sigma-1}$ which is an image of a projection
 $e$ of $A$.
The existence of such a projection $e$ follows since $A$
is AF.
Let $J(e)$ denote the ideal of $A$ generated by $e$ and note
that $A=J(e)+I_{\sigma-1}=J(e)+(1-e)I_{\sigma-1}(1-e)$. Since
$(1-e)I_{\sigma-1}(1-e)$ is an ideal of $(1-e)A(1-e)$ and is
generated by an increasing sequence $(p_n)$ of projections  it
follows that the sequence $J(e+p_n)$ is increasing and $\bigcup_n
J(e+p_n)$ is dense in $A$.
Set $e_n=e+p_n$.

Let $\alpha$ be a flow on $A$.
Then there is an $\alpha$-cocycle $u$
in $A$ (or $A+\C1$ if $A$ is not unital) such that
$\Ad\,u_t\alpha_t(e_n)=e_n$.
To prove that $\alpha$ is
quasi-diagonal on $J(e_n)$ it suffices to show, by Proposition~\ref{pert}, that
$\Ad\,u\alpha$ is quasi-diagonal on $B=e_n Ae_n$.
Note that $B$ has
a composite series $\{J_\mu\}$ with $J_\mu=e_n I_\mu e_n\subset
I_\mu$ and $\mu\in [0,\sigma]$.
Since $B/J_{\sigma-1}\cong \C1$ and
any flow on $J_{\sigma-1}$ is quasi-diagonal one can conclude that
any flow on $B$ is quasi-diagonal.
Thus we conclude that
$\alpha|_{J(e_n)}$ is quasi-diagonal for any $n$, which implies that
$\alpha$ is quasi-diagonal.
\end{pf}

\begin{remark}
If $A$ is a type I AF algebra then any flow on $A$ is an approximate
AF flow $($or a cocycle perturbation of an AF flow $\cite{K05})$. The
proof is quite similar to the above but using Corollary $1.6$ of
$\cite{K06}$ in place of Proposition~$\ref{Her}$. Thus the above lemma
also follows from this fact and Proposition~$\ref{AF}$.
\end{remark}

\begin{remark}
Let $A$ be a type I \cstar\ and  $\alpha$  a flow on $A$ which
fixes each ideal.
Then $\alpha$ is universally weakly inner and is
approximately inner. $($N.B.\ $A$ need not be quasi-diagonal.$)$
This follows from {\rm  \cite{K82,BE85}}.
\end{remark}

\noindent{\bf Problem} {\em Is a quasi-diagonal flow on an AF algebra
approximately inner?}

\medskip

Summarizing Example \ref{GICAR} gives the following.

\begin{prop}\label{TI}
Let $A$ be an AF algebra and suppose that there is a sequence
$\{I_n\}$ of ideals of $A$ such that $\bigcap_nI_n=\{0\}$ and
$A/I_n$ is of type I for all $n$. Then any flow on $A$ is
quasi-diagonal.
\end{prop}

\begin{prop}
Let $A$ be a unital  AF algebra and let $\alpha$ be a flow on $A$.
Suppose that there is a covariant irreducible representation
$(\pi,U)$ such that $(\pi(A),U)$ is quasi-diagonal.
Then there is an
increasing sequence $(E_n)$ of finite-rank projections on $\Hil_\pi$
and an $\alpha$-cocyle $u$ in $A$ such that $[E_n, \pi(u_t)U_t]=0$
and $\bigcup_nB_n$ is dense in $A$ where
$$B_n=\{x\in A;\
[E_k,\pi(x)]=0,\ k\geq n\}.
$$
\end{prop}
\begin{pf}
Let $(A_n)$ be an increasing sequence of finite-dimensional
C$^*$-subalgebras of $A$ such that the union is dense in $A$.
We omit $\pi$ in the arguments below.

Given a finite-rank projection $E$ on $\Hil$ and $\eps>0$ there
is a finite-rank projection $F$ on $\Hil$ such that $A_1E\Hil\subset
F\Hil$, $\|[F,x]\|\leq \eps \|x\|$ for $x\in A_1$  and
$\|[F,U_t]\|<\eps$ for $t\in [-1,1]$.
Since $A_1$ is finite-dimensional the average of $vFv^*$ over  $v$ in the unitary
group of $A_1$ is in a small vicinity of $F$ (depending on $\eps$
and $\dim A_1$)  which yields a projection $F'$ with $\|F-F'\|$
small by functional calculus.
We note that $E\leq F'$, $F'\in A_1'$,
and $\|[F',U_t]\|<\eps'+2\|F'-F\|$ for $t\in [-1,1]$. Since
$\|F-F'\|$ can be made arbitrarily small  we now suppose that the
finite-rank projection $F$ satisfies $E\leq F$, $F\in A_1'$
and $\|[F,U_t]\|<\eps$ for $ t\in [-1,1]$. By
Proposition~\ref{smear} one obtains a projection $F'$ in a small vicinity of $F$
such that $\|[F',H]\|$ is small  where $H$ is the self-adjoint
operator with $U_t=e^{itH}$.
Using the irreducibility of $A$ and Kadison's transitivity theorem we
obtain  an $h=h^*\in A$ such
that $[h,F']=[H,F']$ and $\|h\|=\|[H,F']\|$.
We define $V$ to be the
unitary part of the polar decomposition of $X=F'F+(1-F')(1-F)$.
Since $\|X-1\|\leq 2\|F-F'\|$ and $XF=F'X$ we conclude that
$\|V-1\|$ is small and that $VFV^*=F'$.
A second application of  Kadison's transitivity theorem
gives  a unitary $v\in A$ such that $vF'=VF'$ and
$\|v-1\|\leq \|V-1\|$.
Note that $ve^{it(H-h)}v^*$ commutes with $F$.
Let $u_t=ve^{it(H-h)}v^*e^{-itH}\in A$.
This is an  $\alpha$-cocycle satisfying  $u_tU_tF=Fu_tU_t$ and
$\|u_t-1\|$ is small for $t\in [-1,1]$ because $\|h\|$ and $\|v-1\|$
are small.
Note that if $E$ satisfies that $[E,U_t]=0$ then we may
suppose that $u_tE=E$.

We  apply the foregoing  procedure  repeatedly and each time  make
a perturbation by selecting a cocycle.

Let $(\xi_n)$ be an orthonormal basis for $\Hil$.
We construct a sequence $(E_n)$ of finite-rank projections and a sequence
$(u^{(n)})$ of cocycles such that $E_n\xi_n=\xi_n$, $E_{n-1}\leq
E_n$, $E_n\in A_n'$, $u_t^{(n)}E_{n-1}=E_{n-1}$,
$\|u_t^{(n)}-1\|<2^{-n}$ for $t\in [-1,1]$, and $u^{(n)}$ is
$\alpha^{(n-1)}$-cocycle, where $\alpha^{(0)}=\alpha$ and
$\alpha^{(n)}=\Ad\,u^{(n)}\alpha^{(n-1)}$.
Then $\alpha^{(n)}$ converges to a flow $\alpha^{(\infty)}$ which is a cocycle
perturbation of $\alpha$ with the cocycle $u_t$ obtained as the
limit of $u_t^{(n)}u_t^{(n-1)}\cdots u_t^{(1)}$.
Since $B_n\supset A_n$, the union $\bigcup_nB_n$ is dense in $A$.
The other requirements follow easily.
\end{pf}

Let $\alpha$ be  a flow on a unital simple AF algebra $A$ and
consider the following conditions.
\begin{enumerate}
\item There is a covariant irreducible representation $(\pi,U)$ of
$A$ such that $(\pi(A),U)$ is quasi-diagonal
\item There exists an $\alpha$-cocycle $u$ and an increasing sequence
$(B_n)$ of residually finite-dimensional C$^*$-subalgebras of $A$
such that $\bigcup_nB_n$ is dense in $A$ and each $B_n$ and its
ideals are left invariant under $\Ad\,u\alpha$.
\item $\alpha$ is quasi-diagonal.
\end{enumerate}
The above proposition shows that (1) implies (2).
It is immediate that (2) implies (3) since $\Ad u\alpha|_{B_n}$ is quasi-diagonal.
But to show that (3) implies (1) we would need to extract more
information on $(A,\alpha)$ in addition to the conclusion of
Proposition~\ref{KMS}.

\section{Voiculescu's Weyl-von Neumann theorem}

Our aim is to prove a version of Voiculescu's non-commutative
Weyl-von Neumann theorem \cite{Voi76}. A new feature is that we deal
with a \cstar\  together with a derivation implemented by an unbounded self-adjoint operator.

\begin{theo}\label{WV}
Let $\alpha$ be a flow on a separable \cstar\ $A$ and $(\pi,U)$
$($resp.\  $(\rho,V)$$)$  a covariant representation of $(A,\alpha)$
on a separable Hilbert space such that the range of $\rho\times V$
does not contain a non-zero compact operator. If $\Ker(\rho\times
V)\subset \Ker(\pi\times U)$, then there is a sequence
$W_1,W_2,\ldots$ of isometries from $\Hil_\pi$ into $\Hil_\rho$ such
that $\pi(x)-W_n^*\rho(x)W_n\in \K(\Hil_\pi)$ and
$\|\pi(x)-W_n^*\rho(x)W_n\| \to0$ for $x\in A$. In addition
$H_U-W_n^*H_V W_n\in \K(\Hil_\pi)$ and $\|H_U-W_n^*H_V W_n\| \to0$
where $H_U$ $($resp.\  $H_V$$)$  is the self-adjoint generator of
$U$ $($resp.\ $V$$)$. Furthermore if $\Ker(\rho\times
V)=\Ker(\pi\times U)$ and the range of $\pi\times U$ does not
contain a non-zero compact operator then the $W_n$ can be assumed to
be unitary.
\end{theo}

This theorem immediately allows a proof of Theorem~\ref{Voi}.

\smallskip

\begin{pff}$\;${\em of Theorem~\ref{Voi}.} $\,$
%{\em Proof of Theorem \ref{Voi}}
Suppose that $\alpha$ is a quasi-diagonal flow on a \cstar\ $A$.
Thus there is a  faithful covariant representation $\pi$ of $A$ and
a unitary flow $U$ on the same space satisfying
$\Ad\,U_t\pi(x)=\pi\alpha_t(x)$ such that $(\pi(A),U)$ is
quasi-diagonal. By taking the direct sum of $(\pi,\chi_pU)$ with $p$
rational,  if necessary, we may further suppose that $\pi\times U$
is a faithful representation of $A\times_\alpha\R$ and the range of
$\pi\times U$ does not contain a non-zero compact operator. Let
$(\rho,V)$ be another pair such that $\Ker(\rho\times V)=\{0\}$ and
$\Ran(\rho\times V)\cap \K(\Hil_\rho)=\{0\}$. Let $\F$ be a finite
subset of $A$, $\omega$ a finite subset of $\Hil_\rho$  and
$\eps>0$. Let $B$ be the $\alpha$-invariant \cstarsub\ of $A$
generated by $\F$. Since $B$ is a separable \cstar\   we can apply
Theorem~\ref{WV}  to $(\pi|_B,U)$ and $(\rho|_B,V)$. Since
$\Ker(\pi|_B\times U)=\{0\}=\Ker(\rho|_B\times V)$ and $\pi|_B\times
U$ and $\rho|_B\times V$ satisfy the range condition there is a
unitary $W$ from $\Hil_\rho$ onto $\Hil_\pi$ such that
$\|\rho(x)-W^*\pi(x)W\|\leq(\eps/4)\|x\|$ for $x\in \F$ and
$\|H_V-W^*H_U W\|<\eps/2$. Let $E$ be a finite rank projection on
$\Hil_\pi$ such that $\|[E,\pi(x)]\|\leq (\eps/2)\|x\|$ for $x\in
\F$, $\|(1-E)W\xi\|\leq \eps\|\xi\|$ for $\xi\in\omega$ and
$\|[E,H_U]\|<\eps/2$. We set $F=W^*EW$, which is a finite-rank
projection on $\Hil_\rho$. Then $\|[F,\rho(x)]\|\leq\eps\|x\|$ for
$x\in\F$, $\|(1-F)\xi\|=\|(1-E)W\xi\|\leq \eps\|\xi\|$ for
$\xi\in\omega$, and $\|[F,H_V]\|<\eps/2+\|[E,H_U]\|<\eps$. This
completes the proof in the quasi-diagonal case.
A similar proof can be given in the pseudo-diagonal case.
\end{pff}

Now we turn to the proof of Theorem \ref{WV}.

\smallskip

Let $\sigma$ be a completely positive map (or CP map) of
$A\times_\alpha\R$ into $B(\Hil)$.
Let $(e_\nu)$ be an approximate
identity in $A\times_\alpha\R$.
Then $\sigma(e_\nu)$ is increasing
and bounded in $B(\Hil)$ and hence converges in the strong operator
topology.
Denote the limit by $I$ and remark that $I$ is the
supremum of $\{\sigma(x);\ x\in A\times_\alpha\R,\ 0\leq x\leq 1\}$
in $\B(\Hil)$.

More generally we extend $\sigma$ to a CP map from the multiplier
algebra $M(A\times_\alpha\R)$ into $B(\Hil)$. For $x\in
M(A\times_\alpha\R)$ one shows that $\sigma(e_\nu x)$ converges in
the weak operator topology since $|\lan \xi,
\sigma((e_\mu-e_\nu)x)\eta\ran|\leq
\lan\xi,\sigma(e_\mu-e_\nu)\xi\ran^{1/2}\lan\eta,
\sigma(x^*(e_\mu-e_\nu)x)\eta\ran^{1/2}$ for $\mu\geq \nu$. We
denote the limit by $\sigma(x)$. It  is also the limit of
$\sigma(e_\nu xe_\nu)$ since $\sigma(e_\nu x(e_\mu-e_\nu))$
converges to zero for all pairs $\mu,\nu$ with $\mu>\nu$ as $\nu
\to\infty$. Note that $\sigma(\lambda(f))$ for $f\in L^1(\R)$ and
$\sigma(\lambda_t)$ for $ t\in \R$ are all well-defined where
$\lambda$ denotes the unitary group implementing $\alpha$ on $A$ and
$\lambda(f)=\int f(t)\lambda_tdt\in A\times_\alpha\R$. From the
definition we have that $\sigma(1)=I$.

We next define the $\alpha$-spectrum of $\sigma$.
Set $J=\{f\in L^1(\R);\ \lambda(f)\geq 0,\ \sigma(\lambda(f))=0\}$.
Then $J$ is a hereditary closed cone in $L^1(\R)_+$ and the
$\alpha$-spectrum, denoted by $\Sp_\alpha(\sigma)$, of $\sigma$ is
defined by $\bigcap\{\Ker \hat{f};\ f\in J\}$. When
$\Sp_\alpha(\sigma)$ is compact and $f\in L^1(\R)_+$ satisfies
$\hat{f}=1$ on $\Sp_\alpha(\sigma)$  then it follows that
$\sigma(\lambda(f))=\sigma(1)$.

Note that $t\mapsto \sigma(\lambda_t)$ need not be continuous (even
if $\sigma$ is a state).
Let $\D$ denote the set of $\xi\in \Hil$
such that $\sigma(\lambda_t-1)\xi/it$  converges strongly as $t \to0$.
If $\D$ is a dense linear subspace  then the operator $H'$ defined
on $\D$ as the limit of $\sigma(\lambda_t-1)/it$ is
symmetric.
If the closure of $H'$, which we will denote by
$\sigma(H)$, is self-adjoint  we will say that $\sigma$ is
$\alpha$-differentiable.
If $\Sp_\alpha(\sigma)$ is compact then
$\sigma$ is $\alpha$-differentiable and $\sigma(H)$ is bounded
(because  if $f\in L^1(\R)$ satisfies $\hat{f}=1$ on
$\Sp_\alpha(\sigma)$ and $\supp(\hat{f})$ is compact then it follows
that $\sigma(\lambda_t)=\sigma(\lambda(f)\lambda_t)$ and $t\mapsto
\lambda(f)\lambda_t$ is differentiable in $t$).
If $\sigma$ is a
homomorphism  then $t\mapsto \sigma(\lambda_t)$ is a unitary flow
and thus $\sigma$ is $\alpha$-differentiable.

From now on we always assume that the \cstar\ $A$ is separable.

\begin{definition}
Let $\sigma\colon A\times_\alpha\R \to B(\Hil)$ and
$\sigma'\colon A\times_\alpha\R \to \B(\Hil')$ be $\alpha$-differentiable
CP maps.
For two bounded $($or unbounded self-adjoint$)$ operators $T$
and $T'$  $($with a common domain$)$ we write  $T\thicksim T'$ if the
difference $T-T'$ is  $($or extends to$)$ a compact operator.

We write $\sigma\thicksim\sigma'$ if there is a unitary $V\colon\Hil \to
\Hil'$ such that $\sigma(x)\thicksim V^*\sigma'(x)V$ for $x\in A\cup
A\times_\alpha\R$ and $\sigma(H)\thicksim V^*\sigma'(H)V$  and
$\sigma\lesssim\sigma'$ if there is an isometry $W\colon \Hil \to\Hil'$
such that $\sigma(x)\thicksim W^*\sigma'(x)W$ for $x\in A\cup
A\times_\alpha\R$ and $\sigma(H)\thicksim W^*\sigma'(H)W$  where
$W^*\D(\sigma'(H))=\D(\sigma(H))$ and $WW^*\D(\sigma'(H))\subset
\D(\sigma'(H))$.

We write $\sigma\thickapprox\sigma'$ if there is a sequence of
unitaries $V_n\colon \Hil \to\Hil'$ such that the $V_n$ satisfy the above
conditions for $V$ and $\|\sigma(x)-V_n^*\sigma'(x)V_n\| \to0$ for
$x\in A\cup A\times_\alpha\R$ and
$\|\sigma(H)-V_n^*\sigma'(H)V_n\| \to0$, and
$\sigma\lessapprox\sigma'$ if there is a sequence of isometries
$W_n\colon \Hil \to\Hil'$ such that the $W_n$ satisfy the above conditions
for $W$, $W_nW_n^* \to0$, $\|\sigma(x)-W_n\sigma'(x)W_n\| \to0$ for
$x\in A\cup A\times_\alpha\R$ and
$\|\sigma(H)-W_n^*\sigma'(H)W_n\| \to0$.
\end{definition}

Note that $\lessapprox$ is transitive: if $\sigma\lessapprox\sigma'$
and $\sigma'\lessapprox\sigma''$ then $\sigma\lessapprox\sigma''$ as
easily follows.

We are now able to establish the following version
of the Theorem~\ref{WV}.

\begin{theo}\label{WV'}
Let $(\pi,U)$ be a covariant representation of $(A,\alpha)$ on a
separable Hilbert space such that the range of $\pi\times U$ does
not contain a non-zero compact operator. If $\rho$ is an
$\alpha$-differentiable CP map of $A\times_\alpha\R$ into
$\B(\Hil_\rho)$ such that $Q\rho$ is a homomorphism with
$\Ker(Q\rho)\subset \Ker(\pi\times U)$, where $Q$ is the quotient
map of $\B(\Hil_\rho)$ onto $\B(\Hil_\rho)/\K(\Hil_\rho)$, then
$\sigma\times U\lessapprox\rho$.
\end{theo}

The following lemma is an adaptation of 3.5.5 of \cite{HR}.

\begin{lem} \label{VW1}
Let $\sigma$ be a homomorphism of $A\times_\alpha\R$ into $B(\Hil)$
where $\Hil$ is a separable Hilbert space.
Then there exists a sequence of CP maps $\sigma_n\colon A\times_\alpha\R \to B(\Hil_n)$ such
that $\Sp_\alpha(\sigma_n)$ is compact, $\dim \Hil_n$ is finite,
$\sigma_n(1)=1$  and
$\sigma\lessapprox\bigoplus_{n=1}^\infty\sigma_n$.
\end{lem}
\begin{pf}
Let $P$ be the spectral measure of the generator  $\sigma(H)$
of $t\mapsto\sigma(\lambda_t)$.
On each spectral subspace $P(n,n+1]\Hil$ we find a
compact operator $K_n$ such that $H_n'=\sigma(H)P(n,n+1]+K_n$ is
diagonal with eigenvalues in $(n,n+1]$ and $\|K_n\|<1/(|n|+1)$.
Then $H'=\sum_nH_n'$ is diagonal and is given by $\sigma(H)+K$  where
$K=\sum_nK_n$ is compact.
Let $(E_n)$ be an approximate unit for
$\K(\Hil)$ consisting of projections such that $[E_n,H']=0$.
Note that each $E_n$ commutes with $P(k,k+1]$, is dominated by $P[-m,m]$
for some $m$ and satisfies $\|[E_n,H]\|=\|[E_n,K]\| \to0$ as
$n \to\infty$.
In the convex hull of $(E_n)$ there is an approximate
unit $(F_n)$ for $\K(\Hil)$ such that $(F_n)$ satisfies
$\|[F_n,\sigma(x)]\| \to0$ for $x\in A\cup A\times_\alpha\R$ in addition to
the conditions on $H$.
Note that $F_n$ is of finite rank.

Let $(\F_k)$ be an increasing sequence of finite subsets of $A\cup
A\times_\alpha\R$ such that $\bigcup_k\F_k$ is dense in $A\cup
A\times_\alpha\R$ and $\eps>0$.
By assuming
$\|[F_n,\sigma(x)]\|$ and $\|[F_n,\sigma(H)]\|$ are sufficiently
small we may suppose that $D_n=(F_n-F_{n-1})^{1/2}$ with $F_0=0$
satisfies
$$
\|D_n\sigma(x)-\sigma(x)D_n\|<\eps 2^{-n},\ \ x\in \F_n,
$$
and
$$ \|D_n\sigma(H)-\sigma(H)D_n\|<\eps 2^{-n}.
$$
The former follows from $\|[F_n-F_{n-1},\sigma(x)]\|\approx0$ as
$D_n$ is just a continuous function of $F_n-F_{n-1}$ which can be
approximated by polynomials.
For the latter, where $\sigma(H)$ is
unbounded in general, we use the fact that both $D_n$ and
$\sigma(H)$ commute with $P(k,k+1]$.
Thus we have
$$[D_n,\sigma(H)]=\sum_k
[D_nP(k,k+1],(\sigma(H)-k)P(k,k+1]],
$$
where each commutator is between elements of norm one or less.
Since $D_nP(k,k+1]=(F_nP(k,k+1]-F_{n-1}P(k,k+1])^{1/2}$ the latter
inequality then follows just as the former.

Introduce the finite-dimensional subspace $\Hil_n$ of $\Hil$ by  $\Hil_n=D_n\Hil$.
Let $\Hil'=\bigoplus_n\Hil_n$ and define a linear map $V$
from $\Hil$ into $\Hil'$  by $ V\xi=\bigoplus_n D_n\xi$.
This is an isometry since  $\sum_nD_n^2=1$.
Let $Q_n$ be the projection onto
$\Hil_n$ in $\Hil$ and let $\sigma_n(x)=Q_n\sigma(x)Q_n$ for $x\in
A\times_\alpha\R$.
Since $\Hil_n$ is finite-dimensional and
$D_n=D_nP[-m,m]$ for some $m$ one has  $\sigma_n(1)=1$ and
$\Sp_\alpha(\sigma_n)$ is compact with
$\sigma_n(H)=Q_n\sigma(H)Q_n=Q_n\sigma(H)P[-m,m]Q_n$.
Define
$\sigma'=\bigoplus_n \sigma_n$.
Then $\sigma'(H)$ is well-defined
and is equal to $\bigoplus_n\sigma_n(H)$.

We will show that $V\D(\sigma(H))\subset \D(\sigma'(H))$ and
$\D(\sigma(H))\supset V^*\D(\sigma'(H))$. If $\xi\in \D(\sigma(H))$
then
$$\sum_n\|\sigma_n(H)D_n\xi\|^2\leq \sum_n
\|[\sigma_n(H),D_n]\xi+D_n\sigma(H)\xi\|^2\leq
2\sum_n\|[\sigma(H),D_n]\xi\|^2+2\|\sigma(H)\xi\|^2,
$$
which implies that $V\xi\in \D(\sigma'(H))$.
If $\eta\in \D(\sigma'(H))$ and  $P_N$ denotes the
projection onto the first $N$ direct summands in $\Hil'$  then
\BE\|\sigma(H)V^*P_N\eta\|
&=&\|\sum_{n=1}^N[\sigma_n(H),D_n]\eta_n+D_n\sigma_n(H)\eta_n\|\\
&\leq&
\|\sum_{n=1}^N[\sigma_n(H),D_n]\eta_n\|+\|\sum_{n=1}^ND_n\sigma_n(H)\eta_n\|
\\ &\leq&\eps\|\eta\|+\|V^*P_N\sigma'(H)\eta\|,
\EE
where we have  used  $\|[\sigma(H),D_n]\|<\eps 2^{-n}$.
From this kind of computation we can conclude that $(\sigma(H)V^*P_N\eta)_N$
is a Cauchy sequence and so $V^*\eta\in \D(\sigma(H))$.

From the two inequalities above it follows that
$\D(\sigma(H))=V^*\D(\sigma'(H))$ and
$V\D(\sigma(H))=VV^*\D(\sigma'(H))\subset \D(\sigma'(H))$. Thus
$\D(V\sigma(H))=\D(\sigma(H))$ and
$\D(\sigma'(H)V)=V^*\D(\sigma'(H))$ are equal and on this common
domain
$$
V\sigma(H)\xi-\sigma'(H)V\xi=(D_n\sigma(H)\xi-\sigma(H)D_n\xi)_n.
$$
Since $\|[D_n,\sigma(H)]\|<\eps 2^{-n}$, the closure of
$V\sigma(H)-\sigma'(H)V$ on $\D(\sigma(H))$ is a compact operator
with norm less than $\eps$.
Since the closure of
$V\sigma(H)V^*-\sigma'(H)VV^*$ on $\D(\sigma'(H))$ is compact  it
follows that the closure of $V\sigma(H)V^*-VV^*\sigma'(H)$ on
$\D(\sigma'(H))$ is also compact.
Thus we can conclude that the
closure of $\sigma'(H)VV^*-VV^*\sigma'(H)$ on $\D(\sigma'(H))$ is
compact.
In the same way one concludes that $V\sigma(x)-\sigma'(x)V$
is a compact operator for $x\in \bigcup\F_k$ and so for all $x\in A$
and all $x \in A\times_\alpha\R$.
Note that $\|V\sigma(x)-\sigma'(x)V\|<\eps$ for $x\in \F_1$.
This concludes the proof of $\sigma\lesssim\sigma'=\bigoplus_n\sigma_n$.

It also follows from the foregoing construction of  $\sigma'$ and $V$  that
one has bounds $\|V\sigma(H)-\sigma'(H)V\|<\eps$ and
$\|V\sigma(x)-\sigma'(x)V\|<\eps$ for $ x\in \F_1$ where $\F_1$ is a
prescribed finite subset of $A\cup A\times_\alpha\R$.
One  can then
obtain a sequence of such $(\sigma'_k,V_k)$ such that
$\|V_k\sigma(H)-\sigma'_k(H)V_k\| \to0$ and
$\|V_k\sigma(x)-\sigma'_k(x)V_k\| \to0$ for all $x\in A\cup
A\times_\alpha\R$.
Since the direct sum $\bigoplus\sigma'_k$ is of
the form $\bigoplus\sigma_n$ described in the statement this
concludes the proof.
\end{pf}

\begin{lem}\label{VW2}
Let $\rho$ be a homomorphism of $A\times_\alpha\R$ into $\B(\Hil)$
such that $\Ran(\rho)\cap \K(\Hil)=\{0\}$ and let $\sigma$ be a CP
map of $A\times_\alpha\R$ into $\B(\C^n)$ such that $\sigma(1)=1$,
$\Sp_\alpha(\sigma)$ is compact, and $\ker \sigma\supset \ker \rho$.
Then it follows that $\sigma\lessapprox\rho$.

More generally let $\rho$ be an $\alpha$-differentiable CP map of
$A\times_\alpha\R$ into $\B(\Hil)$ such that $Q\rho$ is an
isomorphism where $Q$ is the quotient map from $\B(\Hil)$ onto
$\B(\Hil)/\K(\Hil)$.
Then the same conclusion follows.
\end{lem}
\begin{pf}
We may assume that $A$ is unital.
There is a $C^\infty$-function
$f\in L^1(\R)$ such that the support of $\hat{f}$ is compact,
$0\leq\lambda(f)\leq1$ and $\hat{f}=1$ on $\Sp_\alpha(\sigma)$.
Then it follows that $\sigma(\lambda(f))=1$.
We denote by $\rho_n$
the representation $\id\otimes\rho$ of $M_n(A\times_\alpha\R)$ on
$\Hil_n=\C^n\otimes \Hil$.
We define a state $\phi$ on
$M_n(A\times_\alpha\R)$ by
$$
\phi([x_{ij}])=\frac{1}{n}\sum_{i,j=1}^n\lan
\xi_i,\sigma(x_{ij})\xi_j\ran,
$$
where $(\xi_i)$ is the standard orthonormal basis for $\C^n$. Since
$\phi|_{\Ker(\sigma_n)}=0$, $\Ker(\rho_n)\subset \Ker(\sigma_n)$,
$\phi(\sigma_n(1\otimes\lambda(f)))=1$ and $\Ran(\rho_n)\cap
\K(\C^n\otimes\Hil)=\{0\}$ there is a sequence $(\eta_k)$ of vectors
in $\C^n\otimes\Hil$ such that $\lan \eta_k,\rho_n(x)\eta_k\ran \to
\phi(x)$ for $x\in M_n(A\times_\alpha\R)$,
$\lan\eta_k,\rho_n(1\otimes\lambda(f))\eta_k\ran=1$ and $\eta_k$
converges to zero weakly. Since $\phi(e_{ij}\otimes
\lambda(f))=\delta_{ij}/n$ we may assume
$$\lan\eta_k,\rho_n(e_{ij}\otimes
\lambda(f))\eta_k\ran=\delta_{ij}/n.
$$
If $\eta_k=(\eta_{k1},\ldots,\eta_{kn})\in
\C^n\otimes\Hil\cong\Hil\oplus\Hil\oplus\cdots\oplus\Hil$ then  we define
an isometry $V_k\colon\C^n \to \Hil$ by $V_k\xi_i=\sqrt{n}\eta_{ki}$.
Then
one can conclude that $\|\sigma(x)-V_k^*\rho(x)V_k\| \to0$ for $x\in
A\times_\alpha\R$. (See the proof of 3.6.7 of \cite{HR} for details.)
Since $\rho(\lambda(f))V_k=V_k$, $\sigma(H)=i\lambda(f')$ and
$\rho(H\lambda(f))=i\rho(\lambda(f'))$  it also follows that
$\|\sigma(x)-V_k\rho(x)V_k\| \to0$ for $x\in A$ and
$\|\sigma(H)-V_k\rho(H)V_k\| \to0$.

The condition required for the existence of the foregoing $(\eta_k)$
is precisely the content of the additional statement.
Let $B$ be the \cstar\
generated by the range of $\rho$ and the compact operators.
Then $\sigma Q$ is a CP map of $B$ into $\B(\C^n)$ vanishing on
$\K(\Hil)$.
Then the above arguments show that $\sigma
Q\lessapprox\id_B$ for maps from $B$.
Composing with the homomorphism $\rho$ of $A$ into $B$ one
arrives at  the conclusion.
\end{pf}

The following is an adaptation of 3.5.2 of \cite{HR}.

\begin{lem}\label{VW3}
Let $\rho$ be an $\alpha$-differentiable CP map of
$A\times_\alpha\R$ into $\B(\Hil)$ such that $Q\rho$ is a
homomorphism  where $Q$ is the quotient map of $\B(\Hil)$ onto
$\B(\Hil)/\K(\Hil)$. Let $\sigma_n\colon A\times_\alpha\R
\to\B(\Hil_n)$ be a sequence of CP maps such that
$\dim\Hil_n<\infty$, $\sigma_n(1)=1$ and $\Sp_\alpha(\sigma_n)$ is
compact. If $\sigma_n\lessapprox\rho$ for all $n$ then
$\sigma\equiv\bigoplus_n\sigma_n\lessapprox\rho$.
\end{lem}
\begin{pf}
For each $n$ there is a $C^\infty$-function $f_n\in L^1(\R)$ such
that the support of $\hat{f}$ is compact, $0\leq \lambda(f_n)\leq 1$
and $\hat{f}_n=1$ on $\Sp_\alpha(\sigma_n)$. By
$\sigma_n(\lambda(f_n))=1$ and the assumption
$\sigma_n\lessapprox\rho$ there is a sequence $(V_k)$ of isometries
of $\Hil_n$ into $\Hil$ such that $\rho(\lambda(f_n))V_k=V_k$,
$V_k\eta$ converges to zero weakly for $\eta\in \Hil_n$ and
$\|\sigma_n(x)-V_k^*\rho(x)V_k\| \to0$ for $x\in A\cup
A\times_\alpha\R$ and $x=H$.

Let $(\F_n)$ be an increasing sequence of finite subsets of $A\cup
A\times_\alpha\R$ such that $\bigcup_n\F_n$ is dense in $A\cup
A\times_\alpha\R$ and let $\eps>0$.
We construct inductively  a
sequence $V_n\colon \Hil_n \to \Hil$ of isometries such that
$\rho(\lambda(f_n))V_n=V_n$,
$\|\sigma_n(H)-V_n^*\rho(H)V_n\|<\eps 2^{-n}$ and
$\|\sigma_n(x)-V_n^*\rho(x)V_n\|<\eps 2^{-n}$ for all $x\in
\F_n$, and moreover $V_n\Hil_n$ with $n>1$ is orthogonal to the
finite-dimensional subspace spanned by $V_m\xi, \rho(H)V_m\xi,
\rho(x)V_m\xi, \rho(x^*)V_m\xi$ with $\xi\in \Hil_m,\ x\in \F_n$ and
$m<n$.
(This last condition may require a slight modification of
$V_n$,  retaining $\rho(\lambda(f_n))V_n=V_n$, which will not affect
the condition $\|\sigma_n(H)-V_n^*\rho(H)V_n\|<\eps 2^{-n}$ since
$\rho(H)$ may be replaced by bounded $\rho(H)\rho(\lambda(f_n))$.)
We define $V=\bigoplus_nV_n$, which is an isometry from
$\bigoplus_n\Hil_n$ into $\Hil$. Since for $x\in \F_m$
$$
V^*\rho(x)V=(V_i^*\rho(x)V_j)_{i,j<m} \oplus\bigoplus_{n\geq
m}V_n\rho(x)V_n,
$$
one has
$$
\sigma(x)-V^*\rho(x)V=(\sigma_i(x)\delta_{ij}-V_i^*\rho(x)V_j)_{i,j<m}\oplus
\bigoplus_{n\geq m}(\sigma_n(x)-V_n^*\rho(x)V_n),
$$
where if $m=1$ we ignore the first direct summand, otherwise it is
an operator on $\bigoplus_{i<m}\Hil_i$, which is finite-dimensional.
Then one concludes that the displayed operator is compact and that
$\|\sigma(x)-V^*\rho(x)V\|<\eps$ for $x\in \F_1$. Similarly one
has
$$
\sigma(H)-V^*\rho(H)V=\bigoplus_n(\sigma_n(H)-V_n^*\rho(H)V_n),
$$
which is compact with norm less than $\eps$. \end{pf}

\medskip
\begin{pff}$\;${\em of Theorem~\ref{WV'}.} $\,$
%{\em Proof of Theorem \ref{WV'}}
Let $(\pi,U)$ be a covariant representation of $(A,\alpha)$ and
$\rho$ an $\alpha$-differentiable CP map of $A\times_\alpha\R$ into
$\B(\Hil_\rho)$ as in the theorem.
We may assume that $A$ is unital.

Let $\sigma=\pi\times U$. By Lemma \ref{VW1} we find a sequence
$\sigma_n\colon A\times_\alpha\R \to \B(\Hil_n)$ of CP maps such
that $\dim\Hil_n<\infty$, $\sigma_n(1)=1$, $\Sp_\alpha(\sigma_n)$ is
compact and $\sigma\lessapprox\bigoplus_n\sigma_n$. Since
$\Ker(\sigma_n)\supset \Ker(\sigma)\supset \Ker({\rho})$  Lemma
\ref{VW2} shows that $\sigma_n\lessapprox {\rho}$. Since
$\bigoplus_n\sigma_n\lessapprox{\rho}$ by Lemma \ref{VW3}  one
concludes that $\sigma\lessapprox{\rho}$.
\end{pff}

\begin{pff}$\;${\em of Theorem~\ref{WV}.} $\,$
%{\em Proof of Theorem \ref{WV}}
The first part of the theorem is a special case of Theorem
\ref{WV'}.

Let $\bar{\pi}=\pi\times U$ and $\bar{\rho}=\rho\times V$ and
suppose that $\Ker(\bar{\rho})=\Ker(\bar{\pi})$. Let
$\bar{\rho}^{\infty}$ denote the direct sum of infinite copies of
$\bar{\rho}$. Then applying the first part of the theorem to
$\bar{\rho}^\infty$ and $\bar{\pi}$ we deduce  that
$\bar{\rho}^\infty$ can be approximated by a direct summand
$\bar{\pi}_1$ of $\bar{\pi}$ through a unitary. Here $\bar{\pi}_1$
is obtained as $P\bar{\pi}(\,\cdot\,)P$  where $P$ is a projection
such that $P\D(\bar{\pi}(H))\subset \D(\bar{\pi}(H))$ and
$\|[P,\bar{\pi}(H)]\|$ is small depending on the approximation. Thus
$\bar{\pi}_1$ is an $\alpha$-differentiable unital CP map and this
situation will simply be written as
$\bar{\rho}^\infty\thicksim\bar{\pi}_1$. Writing
$\bar{\pi}_2=(1-P)\bar{\pi}(\,\cdot\,)(1-P)$ one obtains that
$\bar{\pi}\thicksim\bar{\pi}_1\oplus \bar{\pi}_2$ and that
$\bar{\rho}^\infty\oplus \bar{\pi}_2\thicksim \bar{\pi}$. Since
$\bar{\rho}\oplus\bar{\rho}^\infty\thicksim\bar{\rho}^\infty$, one
calculates that $\bar{\rho}\oplus\bar{\pi}\thicksim\bar{\rho}\oplus
\bar{\rho}^\infty\oplus \bar{\pi}_2\thicksim\bar{\rho}^\infty\oplus
\bar{\pi}_2\thicksim\bar{\pi}$. By changing the roles of
$\bar{\rho}$ and $\bar{\pi}$ we conclude that
$\bar{\rho}\thicksim\bar{\pi}$. Since this is true for any degree of
approximation  one obtains the conclusion (see the arguments on page
340 of \cite{Ar} for more details).
\end{pff}

\section{Abstract Characterizations}

Voiculescu \cite{Voi91} gave conditions for \cstars\ to be quasi-diagonal.
By mimicking his  proof  we shall
establish  Theorems~\ref{Ch} and \ref{Ch'}.
\medskip

\begin{pff}$\;${\em of Theorem~\ref{Ch}.} $\,$
%{\em Proof of Theorem \ref{Ch}}
(1)$\Rightarrow$(2).
Suppose  $\alpha$ is quasi-diagonal.
Then there is a faithful representation $\pi$ of $A$ and a unitary flow
$U$ on $\Hil_\pi$ such that $(\pi(A), U)$ is quasi-diagonal.
For any finite subset $\F$ of $A$ and $\eps>0$ there is a finite subset
$\omega$ of unit vectors in $\Hil_\pi$ such that $\max\{|\lan
\xi,\pi(x)\eta\ran|:\ \xi,\eta\in \omega\}\geq(1-\eps/3)\|x\|$
for all $x\in\F$.
Then there is a finite-rank projection $E$ on
$\Hil_\pi$ such that $\|[E,\pi(x)]\|\leq\eps\|x\|$ for $x\in
\F$, $\|(1-E)\xi\|\leq\eps/3$ for $\xi\in \omega$ and
$\|[E,H]\|<\eps$, where $H$ is the self-adjoint generator of $U$.
Set $B=\B(E\Hil_\pi)$, $\beta_t=\Ad\,e^{itEHE}$  and
$\phi(x)=E\pi(x)E$ for $x\in A$.
Then it follows that $\|\phi(x)\|\geq (1-\eps)\|x\|$ for $x\in\F$,
$\|\phi(x)\phi(y)-\phi(xy)\|\leq \eps\|x\|\|y\|$ for $x,y\in\F$
and $\|\beta_t\phi-\phi\alpha_t\|<\eps$ for $t\in [-1,1]$.
The triple $B,\beta,\phi$ satisfies the required conditions.

\medskip

(2)$\Rightarrow$(3).
Let $\F$ be a finite subset of $A$ and $\eps>0$.
Suppose there is a CP map $\phi$ of $A$ into a
finite-dimensional \cstar\ $B$ with a flow $\beta$ such that
$\|\phi\|\leq1$, $\|\phi(x)\|\geq (1-\eps)\|x\|$ and
$\|\phi(xy)-\phi(x)\phi(y)\|\leq\eps\|x\|\|y\|$ for
$x,y\in\F\cup\F^*$ and $\|\beta_t\phi-\phi\alpha_t\|<\eps$ for $t\in
[-1,1]$.
We may assume that $\phi$ is unital. (If $A$ is not unital,
we may assume this by extending $\phi$ as such; if $A$ is unital we
modify $\phi$ using the fact that $\phi(1)$ is close to a
projection.)
By Stinespring's theorem there is a representation
$\pi$ of $A$ and a finite-rank projection $E$ on $\Hil_\pi$ such
that $\phi$ identifies with $E\pi(\,\cdot\,)E$.
It follows that
\BE
\|\phi(x^*x)-\phi(x^*)\phi(x)\|&=& \|E\pi(x^*)(1-E)\phi(x)E\|\\
&=& \|(1-E)\pi(x)E\|^2.
\EE
Since $\|[E,\pi(x)]\|=\max\{\|E\pi(x)(1-E)\|,E\pi(x^*)(1-E)\|\}$  we
obtain (3) except for the condition concerned with the flow.

To prove (3) fully we have to modify $\phi$ and go back to the proof
of Stinespring's theorem.
For a small $\gamma>0$ we replace $\phi$ by
$$
\varphi= \frac{\gamma}{2}\int_{-\infty}^\infty
e^{-\gamma|t|}\beta_{-t}\phi\alpha_tdt.
$$
Then $\varphi(1)=1$, $\beta_{-t}\varphi\alpha_t\leq
e^{\gamma|t|}\varphi$  and
$\|\beta_t\varphi-\varphi\alpha_t\|<\eps$ for $t\in [-1,1]$.
(This method is used in \cite{KR82}.)
Noting that
$\|\beta_{-t}\phi\alpha_t-\phi\|<\eps(1+|t|)$ we have
$\|\varphi-\phi\|<\eps(1+1/\gamma)$, which is an arbitrarily
small constant if $\eps$ is chosen after $\gamma$.
Hence we may assume that $\phi$ satisfies $\beta_{-t}\phi\alpha_t\leq
e^{\eps|t|}\phi$ in addition to the conditions in (2).

The above representation $\pi$ is constructed as follows.
Assuming $B$ acts on a finite-dimensional Hilbert space $\Hil$ we define an
inner product on the algebraic tensor product $A\otimes\Hil$ by
$$
\lan \sum_ix_i\otimes\xi_i,\sum_jy_j\otimes\eta_j\ran=
\sum_{i,j}\lan \xi_i, \phi(x_i^*y_j)\eta_j\ran.
$$
We obtain a Hilbert space $\Hil_\pi$ by the standard process of
dividing $A\otimes \Hil$ out by the null space followed by
completion and then a representation $\pi$ of $A$ on $\Hil_\pi$ from
the multiplication of $A$ on $A$, the first factor of
$A\otimes\Hil$.
Using a flow $V$ on $\Hil$ such that $\beta_t(x)=V_t
xV_t^*$ for $x\in B$ we define a one-parameter group of operators $W_t$
on $\Hil_\pi$ by
$$
W_t\sum_ix_i\otimes \xi_i=\sum_i\alpha_t(x_i)\otimes V_t\xi_i.
$$
Since
$$
\|W_t\sum_ix_i\otimes \xi_i\|^2
= \sum_{i,j}\lan \xi_i, \beta_{-t}\phi\alpha_t(x_i^*y_j)\eta_j\ran,
$$
the group $W$ is well-defined and we obtain  estimates $e^{-\eps
t}1\leq W_t^*W_t\leq e^{\eps t}1$. Denoting the generator of $W$ by
$iK$  one concludes that $\D(K^*)=\D(K)$ and $-\eps\leq-iK^*+iK\leq
\eps 1$. Then the closure  $k$ of $K-K^*$  has norm less than or
equal to $\eps$. Since $W_t\pi(x)W_{-t}=\pi\alpha_t(x)$  we conclude
that $W_t^*W_t\in \pi(A)'$ and so $k\in \pi(A)'$. Set
$U_t=e^{it(K-k/2)}$ for $t\in \R$. The $U$  is a unitary flow such
that $U_t\pi(x)U_t^*=\pi\alpha_t(x)$ for $x\in A$.

We denote by $\Hil'$ the subspace of $\Hil_\pi$ generated by
$1\otimes\xi$ with $\xi\in\Hil$. Let $E$ denote the projection onto
$\Hil'$. Then it follows that $\phi(x)=E\pi(x)E$ for $x\in A$. Since
$W_tEW_{-t}=E$ one obtains $\|U_tEU_t^*-E\|\leq \|k\||t|$. Thus
condition (3) follows.

(3)$\Rightarrow$(1). We  prove this by following the argument  given
in \cite{Voi91}.

Let $(\F_n)$ be an increasing sequence of finite subsets of $A$ with
dense union.
For $(\F_n,n^{-1})$ we choose a covariant
representation $(\pi_n,U^n)$ and a finite-rank projection $E_n$ on
the representation space $\Hil_n$ satisfying the conditions
described in (3).
Let $\Hil=\bigoplus_{n=1}^\infty\Hil_n$,
$\pi=\bigoplus_{n=1}^\infty\pi_n$, $U=\bigoplus_{n=1}^\infty U^n$
and $E=\bigoplus_{n=1}^\infty E_n$.
For each $k\in\N$ let $P'_k$
denote the projection in $\Hil$ onto the first $k$ direct summands
and let $P_k=(1-P'_k)E=E(1-P'_k)$.
If  $\pi_n(H)$ (resp.\ $\pi(H)$) denotes  the self-adjoint generator
of $U^n$ (resp.\  $U$) then
$\pi(H)=\bigoplus_{n=1}^\infty \pi_n(H)$.
Since
$\|[E_n,\pi_n(x)]\|\leq n^{-1}\|x\|$ for $x\in \F_n$ and
$\|[E_n,\pi_n(H)]\|<n^{-1}$  one has $[P_k,\pi(x)]\in \K$ for
$x\in A$ and $[P_k,\pi(H)]\in \K$ where $\K$ denotes the compact
operators on $\Hil$.
If we denote by $\pi\times U$ the covariant
representation of $A\times_\alpha\R$, then it follows that
$\sigma_k=P_k(\pi\times U)P_k$ is an $\alpha$-differentiable CP map
of $A\times_\alpha\R$ into $\B(P_k\Hil)$ such that the composition
$Q\circ\sigma_k$ is an isomorphism and
$$
\sigma_k(H)=(1-P'_k)\cdot\bigoplus_{n=1}^\infty
E_n\pi_n(H)E_n\cdot(1-P'_k),
$$
where $Q$ is the quotient map of $\B(P_k\Hil)$ into
$\B(P_k\Hil)/\K$.

Let $(\rho,V)$ be a covariant representation on a separable Hilbert
space $X$ such that $\rho\times V$ is faithful and  $\Ran(\rho\times
V)\cap \K=\{0\}$. We set $X_\infty=X\oplus X\oplus\cdots$,
$\rho^\infty=\rho\oplus\rho\oplus\cdots$, and $V^\infty=V\oplus
V\oplus\cdots$. Let $G_k$ denote the projection onto the direct sum
of the first $k$ copies of $X$ in $X_\infty$. We will show that
$(\rho^\infty,V^\infty)$ is quasi-diagonal as required in (1).

Fix $k\in\N$.
By Theorem \ref{WV'} applied to
$(\rho^\infty,V^\infty)|_{G_kX_\infty}$ and $\sigma_n=P_n(\pi\times
U)P_n$ we find partial isometries $S_n\colon X_\infty \to \Hil$ such that
$S_n^*S_n=G_k$, $\Ran(S_n)\subset P_n\Hil$ and
$\|S_n\rho^\infty(x)-\pi(x)S_n\| \to0$ for $x\in A$ and
$\|S_n\rho^\infty(H)-\pi(H)S_n\| \to0$.
Similarly from the pair of
$(\pi,U)$ and $(\rho,V)$ we also find isometries $T_n\colon\Hil \to X$
such that $\|T_n\pi(x)-\rho(x)T_n\| \to0$ for $x\in A$ and
$\|T_n\pi(H)-\rho(H)T_n\| \to0$.
We define an isometry $W_n$ of
$\Hil$ into $X_\infty$ by
$$
W_n \xi=S_n^*\xi\oplus
T_n(1-S_nS_n^*)\xi\oplus0\oplus0\oplus\cdots,
$$
which satisfies that $G_kX_\infty\subset W_n\Hil\subset
G_{k+1}X_\infty$ and $\|W_n\pi(x)-\rho^\infty(x)W_n\| \to0$ for $x\in
A$ and $\|W_n\pi(H)-\rho^\infty(H)W_n\| \to0$ as $n \to\infty$. We
note that $G_k\leq W_nP_nW_n^*$.

Now  $F_m=P_m'E$ is a finite-rank projection such that
$[F_m,P_n]=0$ and $F_mP_n=(P_m'-P_n')E \to P_n$ as $m \to\infty$.
Thus choosing $m_n>n$ for each $n$ sufficiently large we may suppose
that $W_nF_{m_n}P_nW_n^*\xi \to\xi$ for each $\xi\in G_kX_\infty$.
Since $F_{m_n}P_n$ commutes with the range of $\sigma_n$,
$W_nF_{m_n}P_nW_n^*$ will serve as the required finite-rank
projection on $X_\infty$ for a large $n$.
\end{pff}
\medskip

\begin{pff}$\;${\em of Theorem~\ref{Ch'}.} $\,$
%{\em Proof of Theorem \ref{Ch'}}
(1)$\Rightarrow$(2).
This is easy.

(2)$\Rightarrow$(3).
Given $(\F,\eps)$  let $\G=\{x,x^*,xy:\ x,y\in\F\}$.
By condition (2) there is a flow $\beta$ on a
finite-dimensional \cstarsub\ $B$ and a CP map $\phi$ of $A$ into
$B$ such that $\|\phi\|\leq1$, $\|\phi(x)\|\geq (1-\eps)\|x\|$
and $\|\phi(x)\phi(y)-\phi(xy)\|\leq \eps\|x\|\|y\|$ for $x,y\in
\G$, and $\|\beta_t\phi(x)-\phi\alpha_t(x)\|\leq \eps\|x\|$ for
$x\in \G$ and $t\in [-1,1]$.
We may assume that $A$ and
$\phi$ are unital.
For  $\gamma=- \log\eps>0$ we replace $\phi$ by
$$
\varphi={\gamma\over 2}\int
^\infty_{-\infty}e^{-\gamma|t|}\beta_{-t}\phi\alpha_tdt.
$$
Since $\|\varphi(x)-\phi(x)\|\leq
(\eps+e^{-\gamma})\|x\|=2\eps \|x\|$ for $x\in \G$, we may
suppose, starting with $\eps/7$ instead of $\eps$, that
$\phi$ satisfies the above properties as well as the covariance
$\beta_t\phi\alpha_t\leq e^{\gamma |t|}\phi$ for $\gamma\approx
-\log\eps$.
We suppose that $B$ acts on a finite-dimensional
Hilbert space $\Hil$ and choose a unitary flow $V$ on $\Hil$ such
that $\beta_t=\Ad\,V_t|_B$.  Then, by Stinespring's construction for
$\phi$ as in the proof of Theorem \ref{Ch}, we obtain a
representation $\pi$ of $A$, a (non-unitary) flow $W$  and a
finite-rank projection $E$ such that $\phi(x)=E\pi(x)E$ for $ x\in A$,
under the identification of $E\Hil$ with $\Hil$,
$W_t\pi(x)W_{-t}=\pi\alpha_t(x)$ for $x\in A$  and $W_tE=EW_tE=V_tE$.
By a perturbation of $W$ we obtain a unitary flow $U$ such that
$\pi\alpha_t(x)=\Ad\,U_t\pi(x),\ x\in A$.
Then we  conclude that $(\pi,U,E,V)$ satisfies the required
properties.

(3)$\Rightarrow$(1).
The proof is similar to the proof of the
corresponding implication in Theorem~\ref{Ch}.
\end{pff}

\medskip

\noindent{\bf Acknowledgement} Part of this work was carried out
whilst the second author was visiting Hokkaido University with
financial support from a JSPS grant of the first author.


\small

\begin{thebibliography}{99}
\bibitem{Ar} W.\  Arveson, Notes on Extensions of C$^*$-algebras, Duke
Math.\ J.\ 44 (1977), 329--355.

\bibitem{BK} B.\ Blackadar, E.\ Kirchberg, Generalized inductive
limits of finite-dimensional \cstars, Math.\ Ann.\ 307 (1997),
343--380.

\bibitem{BEEK}
O. Bratteli, G.A. Elliott, D.E. Evans, A. Kishimoto, Homotopy of a
pair of approximately commuting unitaries in a simple \cstar, J.
Funct. Anal. 160 (1998), 466--523.

\bibitem{BK01} O. Bratteli, A. Kishimoto, AF flows and continuous
symmetries, Rev. Math. Phys. 13 (2001), 1505--1528.

\bibitem{BR1}
O.\  Bratteli  and D.\  W.\  Robinson, {\em Operator algebras and
quantum statistical mechanics}, Vol.\!~1.
\newblock Second edition. Second printing. Springer-Verlag, New York etc., (2002).

\bibitem{BR2}
O.\  Bratteli  and D.\  W.\  Robinson, {\em Operator algebras and
quantum statistical mechanics}, Vol.\!~2.
\newblock Second edition. Second printing. Springer-Verlag, New York etc., (2002).

\bibitem{LB77} L.G. Brown, Stable isomorphisms of hereditary
subalgebras of \cstars, Pacific J. Math. 71 (1977), 335--248.

\bibitem{BE85} L.G. Brown, G.A. Elliott, Universally weakly inner
one-parameter automorphism groups of separable \cstars, II, Math.
Scand. 57 (1985), 281--288.


\bibitem{Br98}
N. P. Brown, AF embeddability of crossed products of AF algebras by
the integers, J. Funct. Anal. 160 (1998), 150--175.

\bibitem{Br04}
N. P. Brown, On quasidiagonal \cstars,  19--64, Adv. Stud. Pure
Math., 38, Math. Soc. Japan, Tokyo, 2004.

\bibitem{EL}
R. Exel, T.A. Loring, Invariants of almost commuting unitaries, J.
Funct. Anal. 95 (1991), 364--376.

\bibitem{HR} N. Higson, J. Roe, {\em Analytic K-Homology}, Oxford
Univ. Press. 2000.

\bibitem{K82} A. Kishimoto, Universally weakly inner one-parameter
automorphism groups of \cstars, Yokohama Math. J. 30 (1982),
141--149.

\bibitem{K00} A. Kishimoto, Locally representable one-parameter
automorphism groups of AF algebras and KMS states, Rep. Math. Phys.
45 (2000), 333--356.

\bibitem{K01} A. Kishimoto, Examples of one-parameter automorphism
groups of UHF algebras, Commun. Math. Phys. 216 (2001), 395--408.

\bibitem{K05} A. Kishimoto, Approximate AF flows, J. evol. equ. 5
(2005), 153--184.

\bibitem{K06} A. Kishimoto, Multiplier cocycles of a flow on a
\cstar , J. Funct. Anal. 235 (2006), 271--296.

\bibitem{KR82} A. Kishimoto, D.W. Robinson, On unbounded derivations
commuting with a compact group of automorphisms, Publ. RIMS, Kyoto
Univ. 18 (1982), 1121--1136.

\bibitem{L97} H. Lin, Almost commuting self-adjoint matrices and
applications, Fields Inst. Commun. 13 (1997), 193--233.

\bibitem{Ped} G.K. Pedersen, {\em \cstars\ and their automorphism
groups}, Academic Press, 1979.

\bibitem{Sak} S. Sakai, {\em Operator algebras in dynamical systems},
Cambridge University Press, 1991.

\bibitem{Voi76} D. Voiculescu, A non-commutative Weyl-von Neumann
theorem, Rev. Roum. Pures et appl. 21 (1976), 97--113.

\bibitem{Voi86} D. Voiculescu, Almost inductive limit automorphisms
and embeddings into AF algebras, Ergod. Th. \& Dynam. Sys. 6 (1986),
475--484.

\bibitem{Voi91} D. Voiculescu, A note on quasi-diagonal
C$^*$-algebras and homotopy, Duke Math. J. 62 (1991), 267--271.

\bibitem{Voi93} D. Voiculescu, Around quasi-diagonal operators, Integr.\
Equat.\  Oper.\  Th.\ 17 (1993), 137--149.

\end{thebibliography}
\end{document}